\documentclass{article}
\usepackage{appendix}
\usepackage{graphicx,float}
\usepackage{amsmath,mathtools,bm}
\usepackage{xcolor}
\usepackage{multirow}
\usepackage{biblatex}
\usepackage{geometry}[margin=1in]
\usepackage{hyperref}
\usepackage{authblk}
\usepackage{tikz}
\usepackage{soul}
\usetikzlibrary{decorations.pathreplacing,positioning,calc,arrows.meta}

\hypersetup{
    colorlinks,
    citecolor=black,
    filecolor=black,
    linkcolor=black,
    urlcolor=black
}

\usepackage{amssymb}
\usepackage{listings}

\title{Physics-informed traditional and generative Mori-Zwanzig--reduced models for multiscale systems}

\author[1,2]{Emily Williams\footnote{Current affiliation: Oden Institute for Computational Engineering and Sciences, University of Texas at Austin}}
\author[1]{Saad Qadeer}
\author[1]{Amanda A. Howard}
\author[1]{Panos Stinis}

\affil[1]{Advanced Computing, Mathematics and Data Division, Pacific Northwest National Laboratory, WA 99354}
\affil[2]{Massachusetts Institute of Technology, MA 02139}

\begin{document}

\maketitle

\begin{abstract}
    Model reduction via the Mori–Zwanzig (MZ) formalism is well understood for coarse-grained representations constructed from standard basis functions such as Fourier functions. We show that the same formalism can be used to build reduced-order models in bases discovered by physics-informed machine learning. The fluctuation-dissipation theorem (FDT) is used to close the resulting non-Markovian evolution equations for the expansion coefficients of the learned basis. We further introduce a complementary, generative closure model. Relative to the Markovian dynamics alone, both the MZ and generative closures substantially improve predictions of the filtered states.
\end{abstract}

\section{Introduction}

For complex multiscale phenomena the fine spatiotemporal resolution required to resolve all scales of motion is prohibitively expensive. Model reduction is a strategy toward computationally efficient simulations. Projection-based reduced-order models (ROMs) aim to derive modified equations for subsets of the necessary degrees of freedom by projecting the equations onto a basis of reduced dimensionality compared to the original solution space \cite{berkooz1993properorthogonal, rathinam2003newlook, chinesta2011shortreview}. Proper orthogonal decomposition (POD) is a data-driven technique that uses the singular value decomposition (SVD) of solution snapshots to obtain a basis to be used in a POD-Galerkin method \cite{kunish2002galerkinproper,chatterjee2000anintroduction}. These methods attempt to capture the essential large-scale dynamics with negligible computational cost relative to full resolution. ROMs can be augmented with a closure model to account for the influence of the truncated modes \cite{aubry_dynamics_1988, ahmed_multifidelity_2023, ahmed_closures_2021, chorin_problem_2005, borggaard_artificial_2011, sirisup_spectral_2004}.

The closure model is part of the ROM which accounts for the subgrid-scale (SGS) dynamics without explicitly resolving them. An accurate SGS model is essential because unresolved dynamics can significantly impact the degrees of freedom that are kept. Such effects can be represented in the ROM through memory terms, which account for the history of the interaction between resolved and unresolved degrees of freedom resulting from the application of the Mori-Zwanzig (MZ) formalism from nonequilibrium statistical mechanics \cite{givon_extracting_2004, mori_transport_1965, zwanzig_memory_1961}. Traditional methods assume that the energy production and dissipation of the small scales are in equilibrium \cite{smagorinsky_general_1963, vreman_eddyviscosity_2004}. These deterministic SGS models, based on this equilibrium assumption, can be overly dissipative, leading to saturation of important dynamical features \cite{sagaut_large_2006, meneveau_scaleinvariance_2000, berselli_mathematics_2006}. In some cases,  traditional eddy-viscosity-based SGS models can be unstable in numerical simulations \cite{givon_extracting_2004, berselli_mathematics_2006, germano_turbulence_1992, lilly_proposed_1992}. These deficiencies arise in part from an inability to capture the backscatter of energy, especially in regions of chaotic and turbulent flow \cite{leith_stochastic_1990, mason_stochastic_1992}.

Interest in developing data-driven and generative stochastic closure approaches has continued to grow, motivated by the need for efficient climate and weather modeling \cite{gottwald_datadriven_2016, perezhogin_generative_2023, arbabi_generative_2022, charalampopoulos_machinelearning_2022}, and developing generative models of underlying dynamics from observations is a rapidly expanding field \cite{volkmann_scalable_2024, finn_generative_2024, ho_denoising_2020, song_score_2020, jacobsen_cocogen_2023}. There is ongoing work toward developing more accurate MZ-based ROMs using ideas from stochastic modeling \cite{price_optimal_2021, palmer_stochastic_2019, cruzeiro_stochastic_2020}. Moreover, stochastic closure approaches can more accurately represent the interactions of the resolved and unresolved scales by accounting for the stochastic effects associated with eliminated degrees of freedom \cite{parish_coarsegraining_2025, berner_stochastic_2017, freitas_importance_2026}. Statistical and machine-learning-based closure strategies for reduced-order modeling in simulating multiscale or chaotic systems remain an open area of research, particularly toward representing unresolved dynamics accurately while also maintaining stability and computational efficiency at fixed cost \cite{ahmed_closures_2021, zhou_turbulence_2021, stinis_renormalized_2015, chorin_discrete_2015, li_incorporation_2015, sanderse_scientific_2024, keyes_multiphysics_2013,wei2026mamba}.

Spectral methods constitute a valuable approach in scientific computing applications toward approximating partial differential equations (PDEs) \cite{tadmor_review_2012, bernardi_spectral_1997, boyd_chebyshev_2000, hesthaven_spectral_2007, meuris_machine-learning-based_2023}. The applicability and effectiveness of spectral methods depend crucially on the choice of basis functions used to expand the solution of a PDE \cite{meuris_machine-learning-based_2023}. In recent years, physics-informed deep learning has emerged as a contender in providing efficient representations of complex functions for the numerical solution of PDEs \cite{alber_integrating_2019, meuris_machine-learning-based_2023, karniadakis_physicsinformed_2021}. These efforts include identifying custom-made candidate basis functions on which to expand the solution of PDEs for evolution using a spectral approach \cite{meuris_machine-learning-based_2023}. Our previous work has demonstrated the benefit of physics-informed operator learning for improving traditional scientific computing methods \cite{williams_what_2024}, by using physics-informed custom basis functions in evolving the Markovian dynamics for model reduction. In this work, we further augment the ROM by applying the MZ formalism and developing a closure model for physics-informed custom basis functions.

\section{Technical approach}

We focus on  the viscous Burgers' equation, written as
\begin{align}
    s_t + s s_x - \nu s_{xx} = 0
\end{align}
on $(x,t) \in (0,2\pi) \times (0,1)$ with periodic boundary conditions. The initial condition is $s_0 = e^{\sin(x)}$ and the viscosity coefficient is $\nu = 0.1$.

\subsection{Application of MZ formalism}

POD--Galerkin-type ROMs provide the prototypical framework for projection-based model reduction. Galerkin projection enforces the governing equations in a weak sense on a reduced subspace. The resulting ROMs evolve a small set of modal coefficients representing the resolved dynamics. Applying the Galerkin projection to a time-dependent PDE results in the weak formulation
\begin{align}
    \langle \tilde{\bm{\phi}}_k, \bm{s}_t + \mathcal{N}[\bm{s}] \rangle = 0
\end{align}
for $1 \leq k \leq p$. We assume we have $p$ orthonormal  basis functions $\{\tilde{\bm{\phi}}_k\}_{k=1}^p$ with respect to a suitable inner product $\langle \cdot, \cdot \rangle$ on $(0,2\pi)$. We can approximate the solution of the PDE using a modal expansion
\begin{align}
    \bm{s}(x,t) \approx \sum_{i=1}^p a_i(t)\,\tilde{\bm{\phi}}_i(x).
\end{align}

Substituting this expansion into the weak formulation and treating the boundary terms appropriately, as described in Appendix \ref{sec:weak-formulation}, yields a system of ODEs for $k = 1,\dots,p$ governing the full-order expansion coefficients
\begin{align}
    \dot{a}_k
    &= -\left\langle \tilde{\bm{\phi}}_k, 
    \mathcal{N}\left[\sum_{i=1}^p a_i(t)\tilde{\bm{\phi}}_i \right]
    \right\rangle \\
    a_k(0) &= \langle \tilde{\bm{\phi}}_k, s_0 \rangle
    = \varepsilon_k
\end{align}

For some cutoff $m < p$, we decompose the full-order coefficients into resolved and unresolved modes as $a = \{\hat{a},\tilde{a}\}$, where $\hat{a} = \{a_k\}_{k=1}^m$ and  $\tilde{a}= \{a_k\}_{k=m+1}^p$ denote the resolved and unresolved modes, respectively. The goal is to develop a reduced model for the dynamics in the resolved set of variables $\hat{a}$ for $m < p$ and compensate for the unaccounted information with memory and noise in addition to the reduced Markovian dynamics.

The most optimal projection operator $P$ in an $L^2$ sense is the expectation conditioned on the resolved variables with respect to a suitable weight function. We take this to be the product
\begin{align}
    \rho(\hat{a},\tilde{a}) = w(\hat{a}) \delta(\tilde{a})
\end{align}
where $w(\hat{a})$ is an arbitrary distribution for the $\hat{a}$'s and $\delta(\tilde{a})$ is the Dirac measure for the $\tilde{a}$. We specifically choose $w$ as the Gaussian measure, defined by 
\begin{align}
    w(\hat{x}) = \prod_{i=1}^m \frac{e^{-\frac{(x_i - \mu_i)^2}{2 \sigma^2_i}}}{\sqrt{2\pi} \sigma_i} 
\end{align}

For viscous Burgers', since all non-constant modes eventually converge to zero, we take $\mu_k = 0$ and $\sigma_k$ to be rapidly decaying, e.g., $\sigma_k = e^{-k}$ for all $k$. We note that this conditional expectation is equivalent to setting the unresolved modes to zero since 
\begin{align}
    P(f)(\hat{a}) = \mathbb{E}[f | \hat{a}] = \frac{\int_{\mathbb{R}^{p-m}} f(\hat{a},\tilde{a}) \rho(\hat{a},\tilde{a}) d\tilde{a}}{\int_{\mathbb{R}^{p-m}} \rho(\hat{a},\tilde{a}) d\tilde{a}} = \frac{\int_{\mathbb{R}^{p-m}} f(\hat{a},\tilde{a}) w(\hat{a}) \delta(\tilde{a}) d\tilde{a}}{\int_{\mathbb{R}^{p-m}} w(\hat{a}) \delta(\tilde{a}) d\tilde{a}}
    = \frac{f(\hat{a},0) w(\hat{a})}{w(\hat{a})} = f(\hat{a},0)
\end{align}

The projection operator $P$ can also be viewed as an orthogonal projection, with respect to the inner product $\langle \cdot, \cdot \rangle_{\mathcal{C}}$, on the subspace of all functions of $\hat{x}$ defined by
\begin{align}
    \langle g_1, g_2 \rangle_{\mathcal{C}} = \int_{\mathbb{R}^p} g_1 (x) g_2(x) \ \rho(x) \,dx.
\end{align}

We can approximate this integral using
\begin{multline}
    \langle g_1, g_2 \rangle_{\mathcal{C}} \approx \sum_{i_1,\dots,i_m=1}^{N_q} g_1(\sigma_1 z_{i_1} + \mu_1, \dots, \sigma_m z_{i_m} + \mu_m,\underbrace{0,...,0}_\text{$p-m$}) \,\cdot \\ g_2 (\sigma_1 z_{i_1} + \mu_1, \dots, \sigma_m z_{i_m} + \mu_m,\underbrace{0,...,0}_\text{$p-m$}) [\omega_{i_1} \dots \omega_{i_m}] \label{CIPApprox}
\end{multline}
where $\{(z_i, \omega_i)\}_{1 \leq i \leq N_q}$ is the Gauss-Hermite quadrature rule for $N_q$ quadrature points.

After projecting a high-dimensional dynamical system onto resolved coordinates with initial conditions $\varepsilon$, the exact reduced dynamics can be expressed as a generalized Langevin equation using the MZ formalism
\begin{align}
    \dot{a}_k(\varepsilon, t) = R_k(a(\varepsilon,t)) + \int_0^t K_k(a(\varepsilon,t-s),s)\,ds + \zeta_k(\varepsilon,t)
\end{align}
where ${R}_k$ is the Markovian term, ${K}_k$ accounts for memory effects, and $\zeta_k$ represents the influence of unresolved initial conditions via the orthogonal dynamics \cite{givon_extracting_2004, mori_transport_1965, zwanzig_memory_1961}.

Applying the projection again eliminates the noise, since it is equivalent to only considering initial conditions of the form $\varepsilon = (\hat{\varepsilon},0)$. The projection operator commutes with nonlinear functions, and so the MZ identity for a resolved variable $a_k(t)$ with initial conditions $a(0) = (\hat{\varepsilon},0)$ is
\begin{align}
    \dot{{a}}_k ((\hat{\varepsilon},0), t) &= R_k(\hat{a}((\hat{\varepsilon},0), t)) + \int_{0}^t K_k(\hat{a}((\hat{\varepsilon},0),t-s),s) \,ds \label{Eq:second_proj}
\end{align}
for $k=1,\dots,m$. The Markovian dynamics are defined by
\begin{align}
    R_k ({\cdot}) = -\left\langle \tilde{\phi}_k, \mathcal{N}\left[\sum_{i=1}^{m} {({\cdot})}_i \tilde{\phi}_i \right] \right\rangle.
\end{align}
The SGS term is therefore the convolution integral representing the memory effects
\begin{align}
    \dot{{a}}_k ((\hat{\varepsilon},0), t) = R_k(\hat{a}((\hat{\varepsilon},0), t)) + \bm{\tau}_k(\hat{a})
\end{align}
where we have defined
\begin{align}\label{eq:tau_integral}
    \bm{\tau}_k(\hat{a}) \equiv \int_{0}^t K_k(\hat{a}((\hat{\varepsilon},0),t-s),s) \,ds.
\end{align}

The technical approach is summarized in the schematic in Fig.~\ref{fig:mz_schematic} and given in more detail in Appendix \ref{sec:weak-formulation}. In this work, we demonstrate stochastic generative modeling approaches toward representing this history-dependent term.

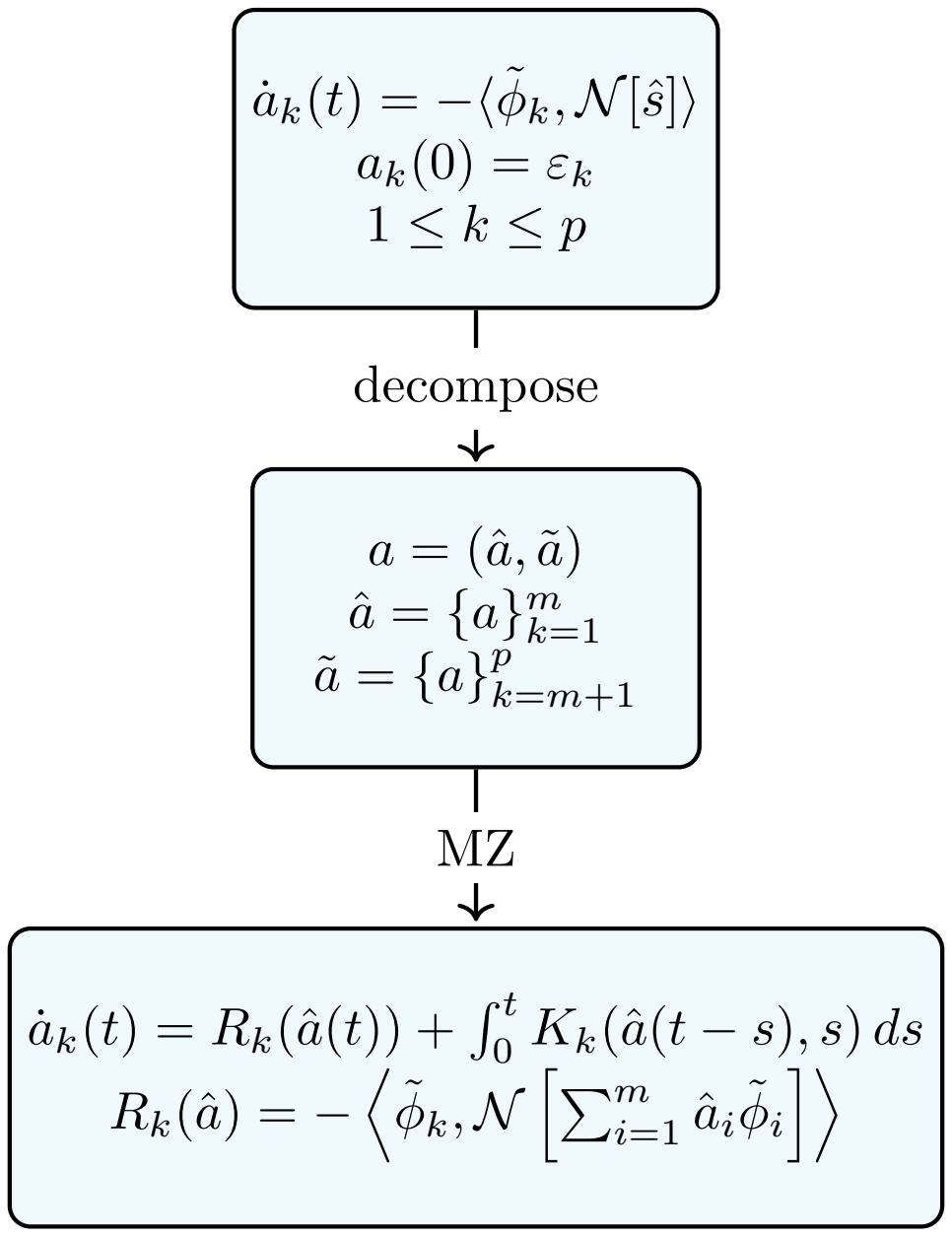
\begin{figure}[H]
\centering

\begin{tikzpicture}[
    box/.style={
        draw,
        rectangle,
        rounded corners,
        align=center,
        minimum width=3cm,
        minimum height=2cm,
        thick,
        fill=cyan!5
    },
    arrow/.style={->, thick, shorten >=1pt},
]

\node[box] (initial) 
{$\dot{a}_k(t) = -\langle \tilde{\phi}_k, \mathcal{N}[\hat{s}] \rangle$ \\
    $a_k (0) = \varepsilon_k$ \\
    $1 \leq k \leq p$};

\node[box, below=3em of initial] (decompose) {
$a = (\hat{a},\tilde{a})$ \\
$\hat{a} = \{a\}_{k=1}^m$ \\
$\tilde{a} = \{a\}_{k=m+1}^p$};

\node[box, below=3em of decompose] (mz) {
    $\dot{{a}}_k(t) = R_k(\hat{a}(t)) + \int_0^t K_k (\hat{a}(t-s),s)\,ds$ \\
    $R_k(\hat{a}) = -\left\langle \tilde{\phi}_k, \mathcal{N}\left[\sum_{i=1}^m \hat{a}_i \tilde{\phi}_i\right]\right\rangle$
};

\draw[arrow] (initial) -- (decompose) node[midway, fill=white] {decompose};
\draw[arrow] (decompose) -- (mz) node[midway, fill=white] {MZ};

\end{tikzpicture}
\caption{Technical approach of applying the MZ formalism to the system of ODEs for the physics-informed expansion coefficients.}
\label{fig:mz_schematic}

\end{figure}

\subsection{Generalized Langevin dynamics}

Our goal is to develop models that can accurately evolve a range of initial conditions. Following \cite{chorin_optimal_2000, chorin_optimal_2002, chorin_problem_2005}, we can expand the memory integrand in terms of kernel functions: if $\{h_l(\hat{x})\}_{l \in I}$ is a basis for functions of the resolved variables and $I$ is the corresponding index set of the chosen basis, then we can write
\begin{align}
    \int_{0}^t K_k(\hat{a}(\varepsilon,t-s),s) \,ds = \sum_{l\in I} \int_0^t \kappa_{lk}(t-s)h_l(\hat{a}(s))\,ds \label{eq:mem_kernels}
\end{align}
for $k = 1,\dots,m$, where $\{\kappa_{lk}(t)\}$ are termed the memory kernels. 

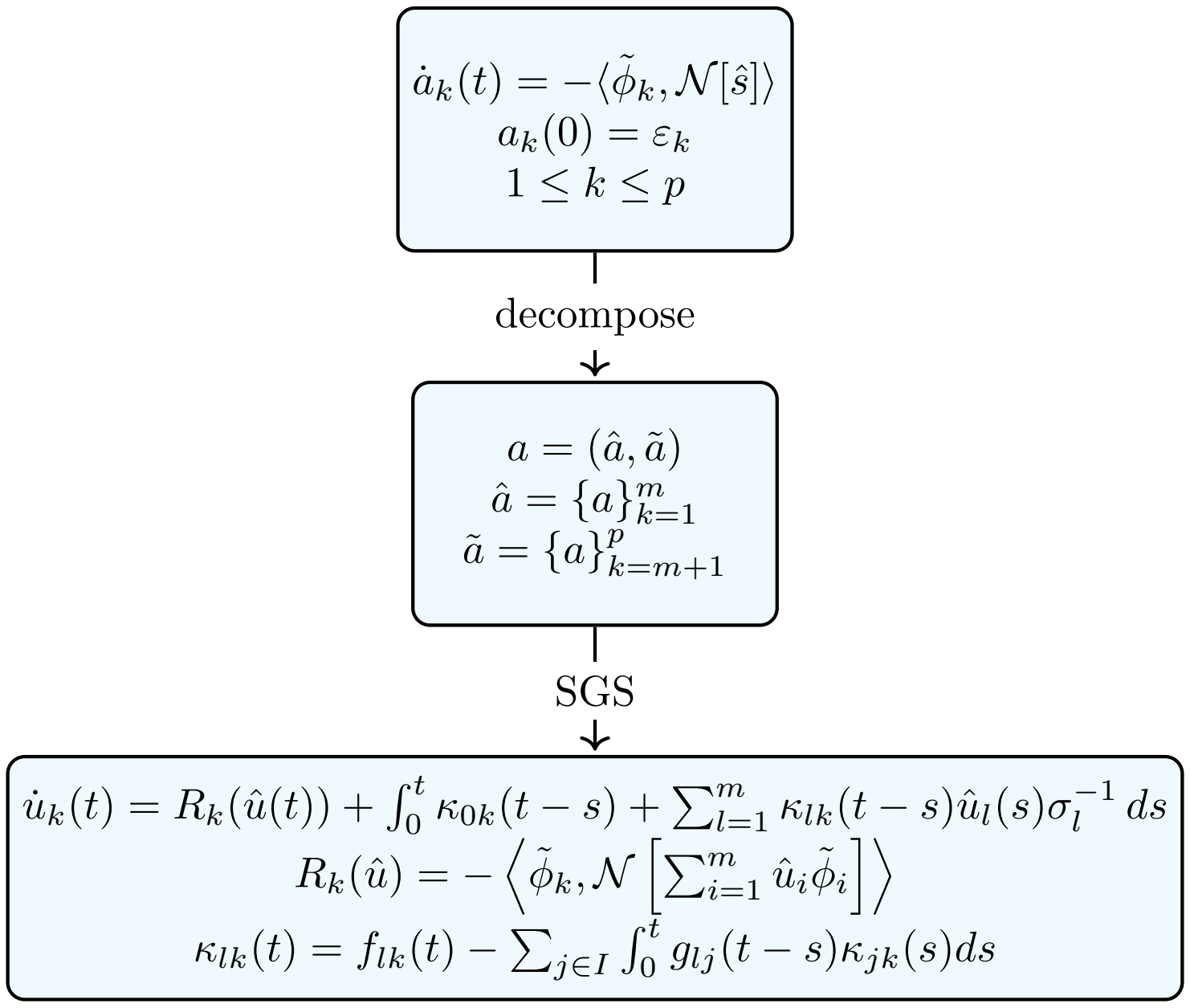
\begin{figure}[H]
\centering

\begin{tikzpicture}[
    box/.style={
        draw,
        rectangle,
        rounded corners,
        align=center,
        minimum width=3cm,
        minimum height=2cm,
        thick,
        fill=cyan!5
    },
    arrow/.style={->, thick, shorten >=1pt},
]

\node[box] (initial) 
{$\dot{a}_k(t) = -\langle \tilde{\phi}_k, \mathcal{N}[\hat{s}] \rangle$ \\
    $a_k (0) = \varepsilon_k$ \\
    $1 \leq k \leq p$};

\node[box, below=3em of initial] (decompose) {
$a = (\hat{a},\tilde{a})$ \\
$\hat{a} = \{a\}_{k=1}^m$ \\
$\tilde{a} = \{a\}_{k=m+1}^p$};

\node[box, below=3em of decompose] (mz) {
    $\dot{{u}}_k(t) = R_k(\hat{u}(t)) + \int_0^t \kappa_{0k}(t-s) + \sum_{l=1}^m \kappa_{lk}(t-s)\hat{u}_l(s) \sigma_l^{-1}\,ds$ \\
    $R_k(\hat{u}) = -\left\langle \tilde{\phi}_k, \mathcal{N}\left[\sum_{i=1}^m \hat{u}_i \tilde{\phi}_i\right]\right\rangle$ \\
    $\kappa_{lk}(t) = f_{lk}(t) - \sum_{j \in I} \int_0^t g_{lj}(t-s) \kappa_{jk}(s) ds$
};

\draw[arrow] (initial) -- (decompose) node[midway, fill=white] {decompose};
\draw[arrow] (decompose) -- (mz) node[midway, fill=white] {SGS};

\end{tikzpicture}
\caption{Technical approach of the FDT implementation for a system of ODEs for the physics-informed expansion coefficients.}
\label{fig:mz_rom_schematic}

\end{figure}

The calculation of these functions is carried out by leveraging the fluctuation-dissipation theorem (FDT) \cite{kubo_fluctuation-dissipation_1966,ceriotti_colored-noise_2010, li_computing_2017}. In practice, we can only use finitely many basis functions so the above expansion is necessarily replaced by an approximation. The technical approach is summarized in the schematic in Fig.~\ref{fig:mz_rom_schematic}.

\section{Traditional MZ-based closure model}

\subsection{Memory kernels}

Following \cite{chorin_optimal_2002}, the memory kernels obey
\begin{align}
    \kappa_{lk}(t) = f_{lk}(t) - \sum_{j \in I} \int_0^t g_{lj}(t-s) \kappa_{jk}(s) ds
\end{align}
where $f_{lk}$ and $g_{lj}$ are auxiliary quantities defined, for $l,j \in I$ and $1\leq k \leq m$, as 
\begin{align}
    f_{lk}(t) = \langle h_l, e^{t\mathcal{L}} \mathcal{L} \mathcal{Q} \mathcal{L} \varepsilon_k \rangle_{\mathcal{C}}, \qquad g_{lj}(t) = \langle h_l, e^{t\mathcal{L}} \mathcal{L} h_j \rangle_{\mathcal{C}} \label{FGdefn}
\end{align}
and where $\mathcal{L} = \sum_{j=1}^p R_j(\varepsilon) \frac{\partial}{\partial \varepsilon_j}$ is the Liouville operator. Denoting by $\{H_l(z)\}_{0 \leq l \leq d}$ the one-dimensional Hermite polynomials corresponding to weight functions $e^{-z^2/2}$, we take $\{h_j(\hat{\varepsilon})\}_{j \in I}$ to be the tensor products of $H_l\left(\frac{\epsilon_k}{\sigma_k}\right)$. In order to establish a proof of concept, we limit ourselves to $d=1$, so we have $|I| = \begin{pmatrix} m+d \\ d \end{pmatrix} = m+1$, where $I$ denote the index set of the multivariate Hermite basis functions.

The inner products in \eqref{FGdefn} are approximated using \eqref{CIPApprox}, with the quadrature points furnishing the initial conditions for generating the full-order solutions. More precisely, these are of the form
\begin{align}
    {a}_k(0) &= \left(\hat{\varepsilon}_k,\tilde{\varepsilon}_k\right) = \left( \left(\sigma_k \bm{z}_k + \mu_k \right)_{k=1}^m, \bm{0}_{k=m+1}^p \right) 
\end{align}

We can compute the physical initial conditions that correspond to these modes via
\begin{align}
    s_0(x) = \sum_{i=1}^m a_i(0) \tilde{\phi}_i(x).
\end{align}

\begin{figure}[H]
    \centering
    \includegraphics[width=0.4\textwidth]{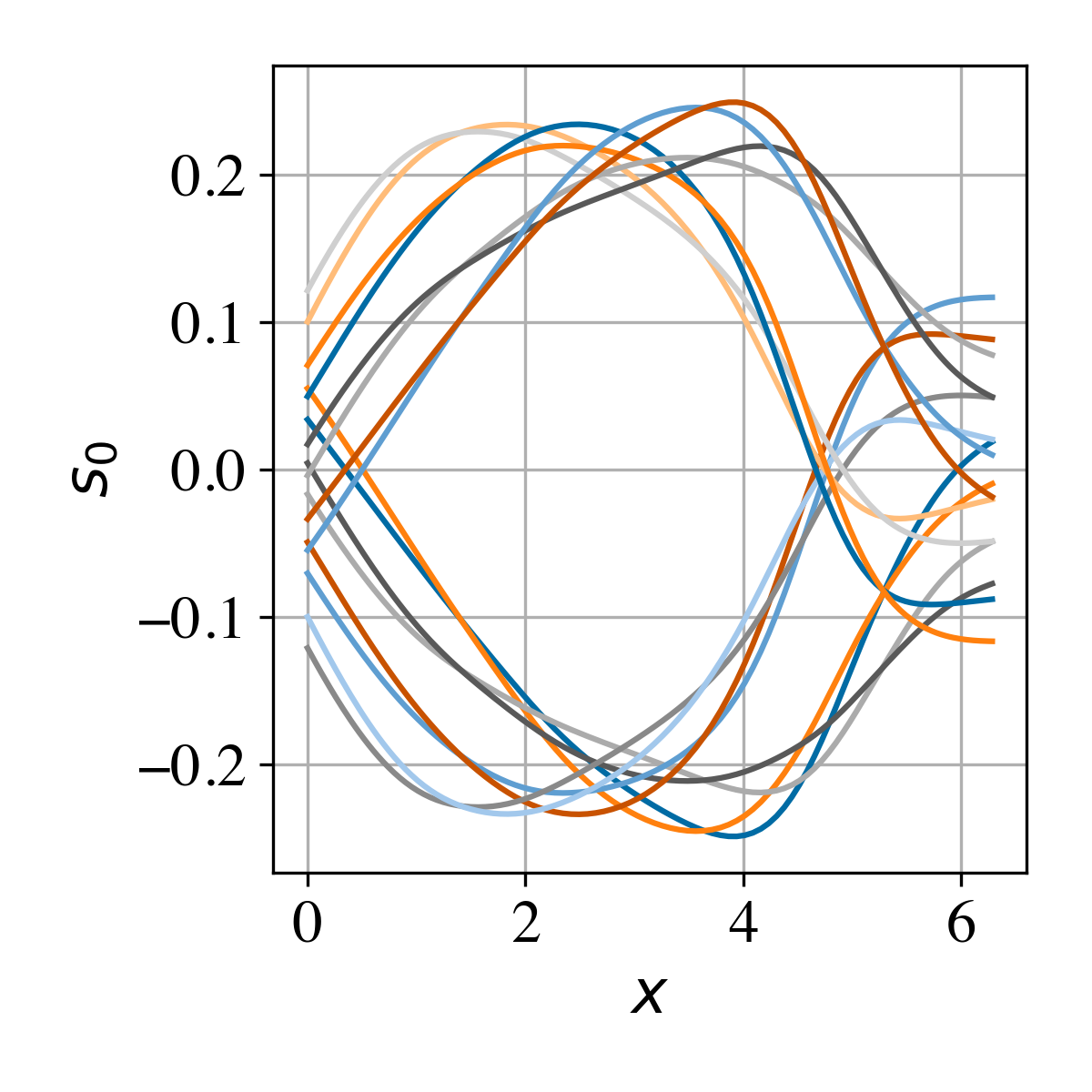}
    \caption{Collection of initial conditions from Gauss-Hermite quadrature space for $N_q = 2$.}
    \label{fig:init-mz}
\end{figure}

The collection of initial conditions for $N_q = 2$ quadrature points is shown in Figure \ref{fig:init-mz}.

\begin{figure}[H]
    \centering
    \includegraphics[width=\textwidth]{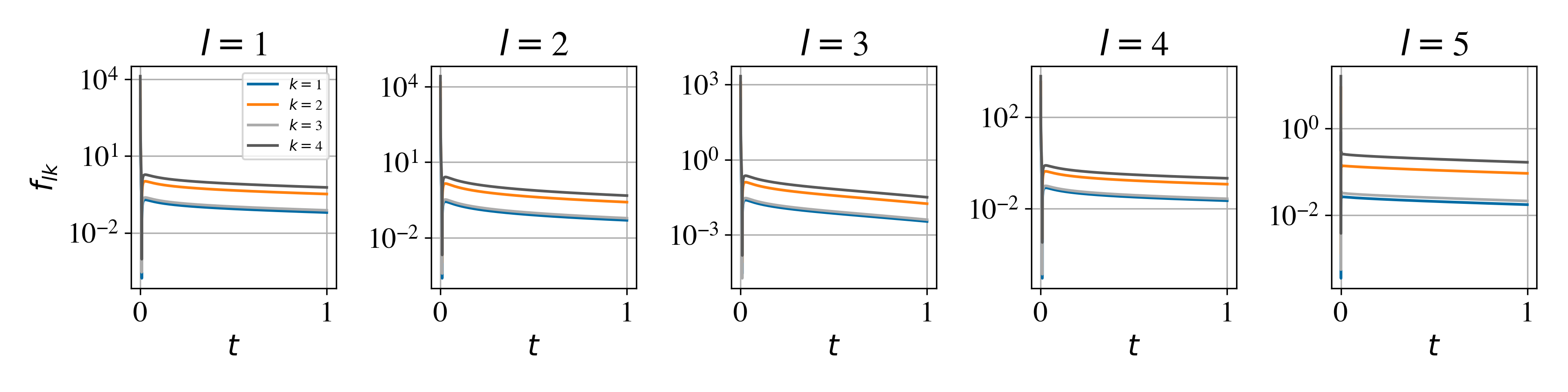}
    
    \includegraphics[width=\textwidth]{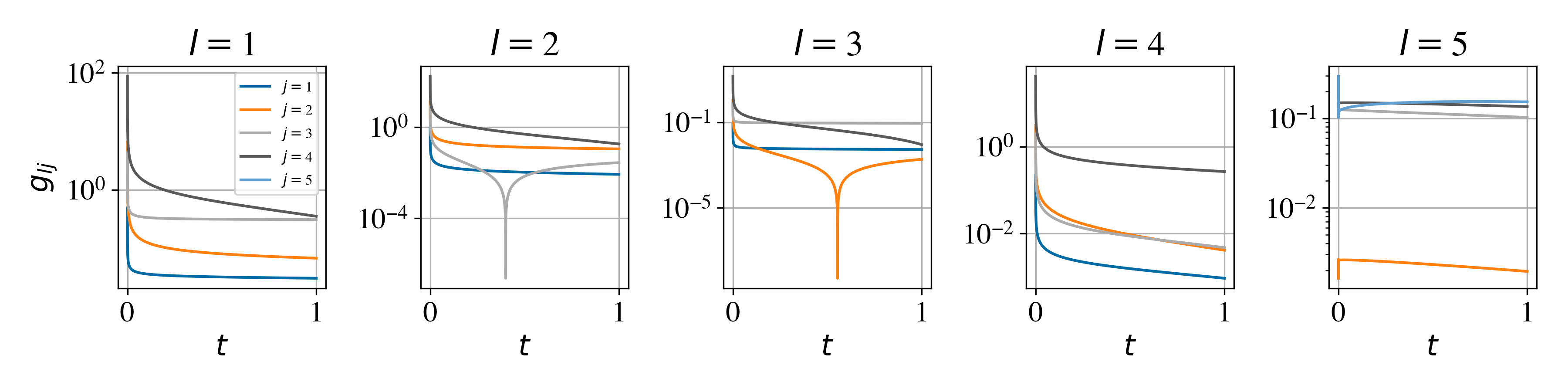}
    \caption{Fluctuations $f_{lk}(t)$ and kernels $g_{lj}(t)$.}
    \label{fig:mz-f-g}
\end{figure}

\begin{figure}[H]
    \centering
    \includegraphics[width=\textwidth]{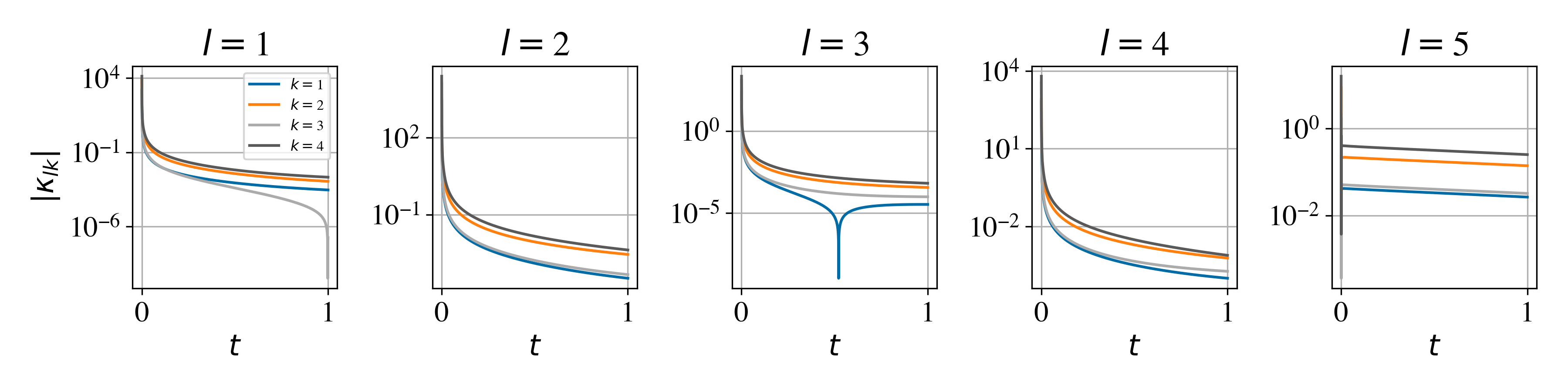}
    \caption{Memory kernels $\kappa_{lk}(t)$.}
    \label{fig:mz-kappa}
\end{figure}

Figure \ref{fig:mz-f-g} plots auxiliary quantities $\bm{f}(t) = (f_{lk}(t))_{l \in I, 1 \leq k \leq m}$ and $\bm{g}(t) = (g_{lj}(t))_{l,j \in I}$. Figure \ref{fig:mz-kappa} plots the memory kernels $\bm{K}(t) = (\kappa_{lk}(t))_{l \in I, 1 \leq k \leq m}$, calculated by discretizing the convolution integral using the trapezoidal rule. See Appendix \ref{sec:mz-kernels} for full derivations. 

\subsection{Reduced model}
For $1 \leq l \leq m$, we model the reduced coefficients $u_k$ and incorporate the memory kernel
\begin{align}
    \dot{u}_k (\varepsilon, t) &= R_k(\hat{u}(\varepsilon, t)) + \int_{0}^t \kappa_{0k}(t-s) + \sum_{l=1}^m \kappa_{lk} (t-s) \frac{\hat{u}_l (s)}{\sigma_l} \,ds
\end{align}
for $k = 1,\dots,m$. The system is integrated forward in time using an explicit trapezoidal rule via Heun's method. 
\begin{align}
    u_k^{n+1} = u_k^n + \frac{\Delta t}{2} (r_1 + r_2)
\end{align}
where
\begin{multline}
    r_1 = R_k(u^n) + \frac{\Delta t}{2} \left[ 2\sum_{q=1}^{n-1} \kappa_{0k}(t_{n-q}) + \kappa_{0k} (t_n) + \kappa_{0k} (t_0) \right]  \\ + \frac{\Delta t}{2} \sum_{l=1}^m \left[ 2 \sum_{q=1}^{n-1} \frac{\hat{u}_l^q}{\sigma_l} \kappa_{lk} (t_{n-q}) + \frac{\hat{u}_l^0}{\sigma_l} \kappa_{lk} (t_n) + \frac{\hat{u}_l^n}{\sigma_l} \kappa_{lk}(t_0) \right]
\end{multline}
\begin{multline}
    r_2 = R_k(u^n + r_1 \Delta t) + \frac{\Delta t}{2} \left[ 2\sum_{q=1}^{n} \kappa_{0k}(t_{n+1-q}) + \kappa_{0k} (t_{n+1}) + \kappa_{0k} (t_0) \right] \\ + \frac{\Delta t}{2} \sum_{l=1}^m \left[ 2 \sum_{q=1}^{n} \frac{\hat{u}_l^q}{\sigma_l} \kappa_{lk} (t_{n+1-q}) + \frac{\hat{u}_l^0}{\sigma_l} \kappa_{lk} (t_{n+1}) + \frac{u_l^n + r_1 \Delta t}{\sigma_l} \kappa_{lk}(t_0) \right]
\end{multline}

Then, the reduced field $\hat{s}$ is reconstructed as
\begin{align}
    \hat{s} = \sum_{i=1}^m {u}_i (t) \tilde{\phi}_i (x)
\end{align}

\begin{figure}[H]
    \centering
    \includegraphics[width=0.8\linewidth]{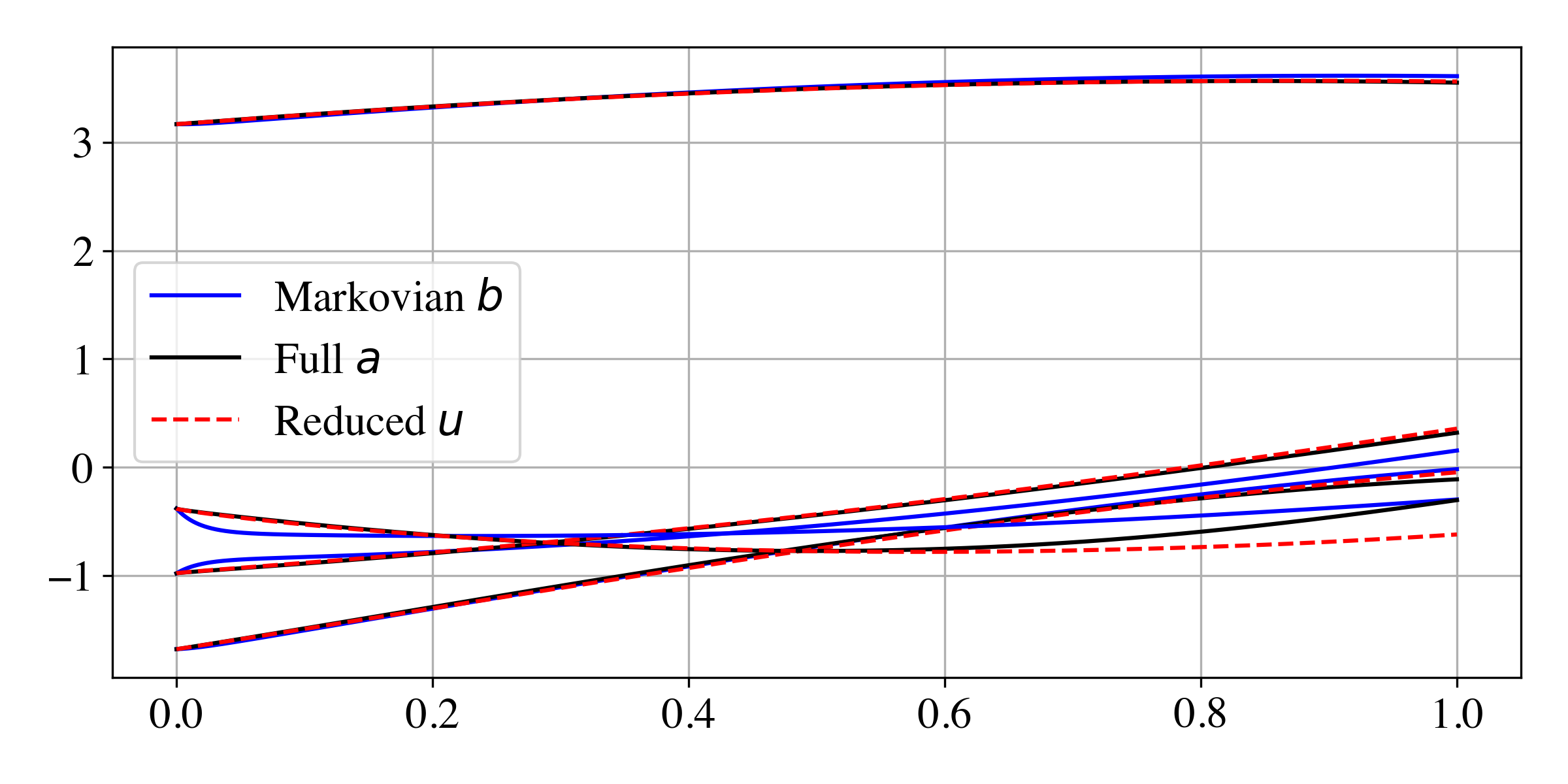}
    \caption{Full-order $a$, reduced $u$, and Markovian $b$ coefficients.}
    \label{rom-gen-coeffs}
\end{figure}

\begin{figure}[H]
    \centering
    \includegraphics[width=\linewidth]{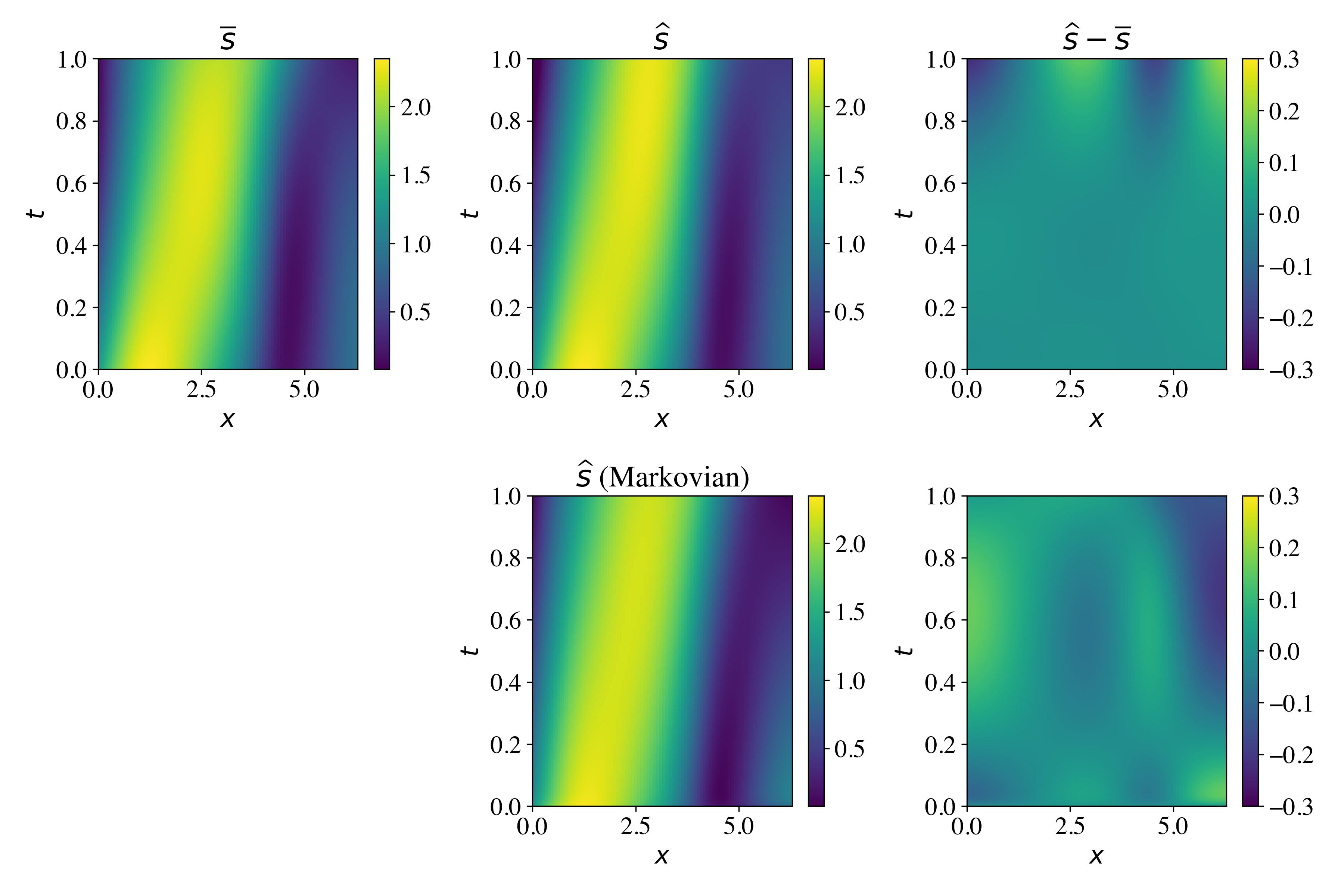}
    \caption{True filtered field $\overline{s}$ compared to Markovian and reduced approximations $\hat{s}$ with absolute error fields for initial condition $s_0 = e^{\sin(x)}$.}
    \label{rom-gen-field}
\end{figure}

Figure \ref{rom-gen-coeffs} plots the evolution of the true full-order coefficients $a$, Markovian coefficients $b$, and reduced coefficients $u$. The reduced dynamics better match the true coefficients than the Markovian dynamics alone. Figure \ref{rom-gen-field} shows the reduced and Markovian fields $\hat{s}$ with absolute errors compared to the true filtered field. We see smaller absolute error magnitudes using the memory kernels compared to the Markovian dynamics alone. We also see increased error towards the end of the temporal domain. These results show that for longer time, we may require higher-order projections (i.e., beyond linear) for representing the memory kernels. We also require a much finer timestep when discretizing the convolution integral, typically smaller by a factor of 100 compared to the full order and Markovian models.

\begin{figure}[H]
    \centering
    \includegraphics[width=0.4\linewidth]{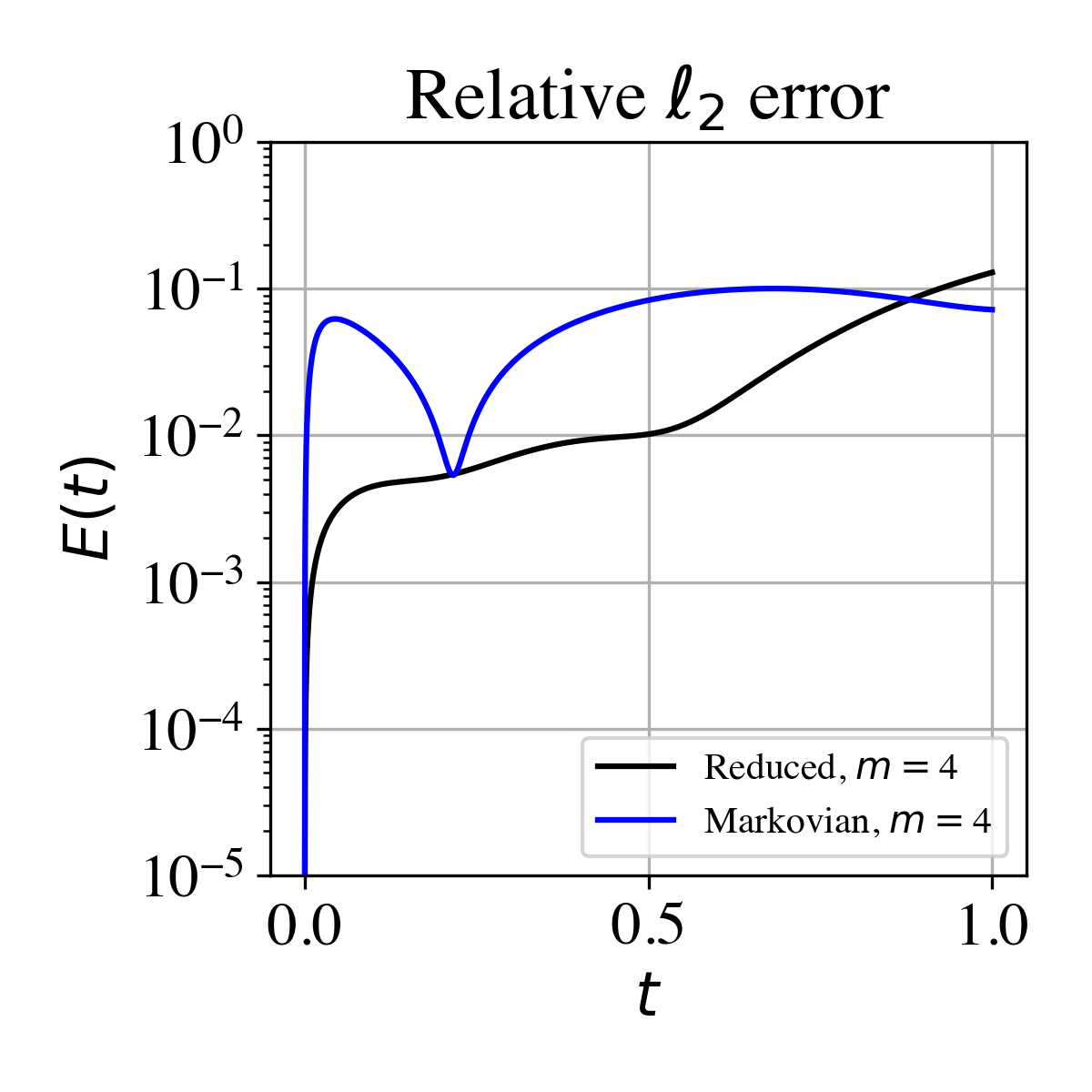}
    \caption{Relative $\ell_2$ error in time.}
    \label{rom-gen-L2}
\end{figure}

Figure \ref{rom-gen-L2} plots the relative $\ell_2$ error of the reduced and Markovian approximations compared to the true filtered field. We see a decreased error across the temporal domain using the memory kernels, except for longer time. Table \ref{tab:viscous-burgers-rom-error} shows the average relative $\ell_2$ error, only 2.9\% for the reduced model compared to almost 7\% for the Markovian model.

\begin{table}[H]
    \centering
    \begin{tabular}{|c|c|}\hline
         & Average relative $\ell_2$ error \\ \hline\hline
    Markovian & 6.6\% \\ \hline
    Reduced & 2.9\% \\ \hline
    \end{tabular}
    \caption{Average errors for independent test sample for the viscous Burgers' equation.}
    \label{tab:viscous-burgers-rom-error}
\end{table}

\section{Generative closure model}

\subsection{Diffusion modeling}

A generative approach samples the SGS term from the data distribution $p_{data}$ along the trajectory in time via learning a denoising diffusion process \cite{flowsanddiffusions2025, song_score_2020, ho_denoising_2020}. We assume access to some initial distribution $p_{init}$ that we can easily sample from, such as the Gaussian $p_{init} = \mathcal{N}(0,I_d)$. The goal of our model is then to transform samples from $\bm{t}_0 \in p_{init}$ into samples from $\bm{t}_1 \sim p_{data}$. We assume access to finite independent samples from $p_{data}$, which serve as a proxy for the true distribution. The simulation of a stochastic differential equation (SDE) can be used for converting a simple distribution $p_{init}$ into a more complex target distribution $p_{data}$. The diffusion model consists of a flow model $u_\gamma^\theta : \mathbb{R}^d \times [0,1] \rightarrow \mathbb{R}^d, (\bm{t},\gamma) \mapsto u_\gamma^\theta(\bm{t})$ with parameters $\theta$ and a fixed diffusion coefficient $\sigma_\gamma : [0,1] \rightarrow [0,\infty), \gamma \mapsto \sigma_\gamma$, such that 
\begin{align}
    \bm{t}_0 &\sim p_{init} = \mathcal{N}(0,I_d) \\
    d\bm{t}_\gamma &= u_\gamma^\theta (\bm{t}_\gamma) d\gamma + \sigma_\gamma dW_\gamma \\
    \bm{\tau}^{\mathrm{SGS}} &= \bm{t}_1 \sim p_{data}
\end{align}

If we set $\sigma_\gamma = 0$, we recover a flow model. The flow training target $u_\gamma^{\mathrm{target}}$ is the marginal vector field.  To construct a training target $u_\gamma^{\mathrm{target}}$ for a flow model, we choose a conditional probability path $p_\gamma(\bm{t}|\bm{\tau})$ that fulfills $p_0(\cdot|\bm{\tau}) = p_{init}$, $p_1(\cdot|\bm{\tau}) = \delta_{\bm{\tau}}$ (Dirac delta distribution). The probability path specifies a gradual interpolation between $p_{init}$ and $p_{data}$. We consider a Gaussian conditional probability path $p_\gamma (\bm{t}|\bm{\tau}) = \mathcal{N}(\bm{t}; \alpha_\gamma \bm{\tau}, \beta_\gamma^2 I_d)$. This is the probability path used by denoising diffusion models \cite{ho_denoising_2020, flowsanddiffusions2025}. Let $\alpha_\gamma : [0,1] \rightarrow \mathbb{R}$ and $\beta_\gamma : [0,1] \rightarrow \mathbb{R}$ be noise schedulers defined as monotonic, continuously differentiable functions satisfying $\alpha_1 = \beta_0 = 1$ and $\alpha_0 = \beta_1 = 0$. Flow matching consists of training a model $u_\gamma^\theta$ via minimizing the conditional flow matching (CFM) loss
\begin{align}
    \mathcal{L}_{\mathrm{CFM}} (\theta) = \mathbb{E}_{\bm{\tau} \sim p_{data}, \gamma \sim \mathrm{Unif}, \bm{t} \sim p_\gamma(\cdot|\bm{\tau})} [ || u_\gamma^\theta (\bm{t}) - u_\gamma^{\mathrm{target}}(\bm{t}|\bm{\tau})||^2 ]
\end{align}
where $u_\gamma^{\mathrm{target}}(\bm{t}|\bm{\tau})$ is the conditional vector field. This can be any vector field such that the corresponding ODE gives the conditional probability path, i.e.,
\begin{align}
    \bm{t}_0 &\sim p_{init} = \mathcal{N}(0,I_d) \\
    d\bm{t}_\gamma &= u_\gamma^{\mathrm{target}}(\bm{t}_\gamma|\bm{\tau}) d\gamma \\
    \bm{t}_\gamma &\sim p_\gamma (\cdot | \bm{\tau})
\end{align}
for $0 \leq \gamma \leq 1$. For a Gaussian probability path, this expression can be found analytically as
\begin{align}
    u_\gamma^{\mathrm{target}}(\bm{t}|\bm{\tau}) = \left(\dot{\alpha}_\gamma - \frac{\dot{\beta}_\gamma}{\beta_\gamma}\alpha_\gamma\right) \bm{\tau} + \frac{\dot{\beta}_\gamma}{\beta_\gamma}\bm{t}
\end{align}

After training, samples are generated by simulating the SDE with stochastic dynamics at each physical timestep for some $\sigma_\gamma$. This form has the same marginal probability path as the corresponding SDE. The score correction ensures that the learned distribution has the correct marginal density evolution toward $p_{data}$ \cite{song_score_2020, ho_denoising_2020, maoutsa_interacting_2020}. We can write the form for the diffusion model as
\begin{align}
    \bm{t}_0 &\sim p_{init} = \mathcal{N}(0,I_d) \\
    d\bm{t}_\gamma &= \left[ u_\gamma^\theta (\bm{t}_\gamma) + \frac{\sigma_\gamma^2}{2} s_\gamma^\theta (\bm{t}_\gamma) \right] d\gamma + \sigma_\gamma dW_\gamma \\
    \bm{\tau}^{\mathrm{SGS}} &= \bm{t}_1 \sim p_{data}
\end{align}
where we train a score network $s_\gamma^\theta$ via denoising score matching (DSM)
\begin{align}
    \mathcal{L}_{\mathrm{DSM}} (\theta) = \mathbb{E}_{\bm{\tau} \sim p_{data}, \gamma \sim \mathrm{Unif}, \bm{t} \sim p_\gamma (\cdot | \bm{\tau})} [ ||s_\gamma^\theta (\bm{t}) - \nabla \log p_\gamma (\bm{t}|\bm{\tau}) ||^2 ]
\end{align}
where the conditional score distribution is
\begin{align}
    \nabla \log p_\gamma(\bm{t}|\bm{\tau}) &= -\frac{\bm{t} - \alpha_\gamma \bm{\tau}}{\beta_\gamma^2}
\end{align}

Further training details and hyperparameters for flow and score networks are in Appendix \ref{sec:training}. At inference time, we can combine $s^\theta_\gamma (\bm{t})$ with $u^\theta_\gamma (\bm{t})$ and sample via simulating the SDE from $\gamma = 0$ to 1. Often, our goal is to generate samples conditioned on some data $\mathcal{S}$. Guided (conditional) generation involves sampling from $\bm{\tau} \sim p_{data}(\cdot | \mathcal{S})$, where $\mathcal{S}$ is a conditioning variable and $p_{data}(\cdot|\mathcal{S})$ is the conditional data distribution. Guided flow matching consists of training a model $u_\gamma^\theta$ via minimizing the classifier-free guidance (CFG) CFM loss
\begin{align}
    \mathcal{L}^{\mathrm{CFG}}_{\mathrm{CFM}} &= \mathbb{E}_\square ||u^\theta_\gamma (\bm{t}|\mathcal{S}) - u^{\mathrm{target}}_\gamma (\bm{t}|\bm{\tau})||^2 \\
    \square &= (\bm{\tau},\mathcal{S}) \sim p_{data}(\bm{\tau},\mathcal{S}), \gamma \sim \mathrm{Unif}, \bm{t} \sim p_\gamma(\cdot|\bm{\tau}), \text{ replace } \mathcal{S} \text{ with } \varnothing \text{ with probability } \eta.
\end{align}

To incorporate score matching, we additionally train a score network $s_\gamma^\theta$ by minimizing the CFG DSM loss
\begin{align}
    \mathcal{L}^{\mathrm{CFG}}_{\mathrm{DSM}} &= \mathbb{E}_\square ||s^\theta_\gamma (\bm{t}|\mathcal{S}) - \nabla \log p_\gamma (\bm{t}|\bm{\tau})||^2 \\
    \square &= (\bm{\tau},\mathcal{S}) \sim p_{data}(\bm{\tau},\mathcal{S}), \gamma \sim \mathrm{Unif}, \bm{t} \sim p_\gamma (\cdot|\bm{\tau}), \text{ replace } \mathcal{S} = \varnothing \text{ with probability } \eta. 
\end{align}

These objectives account for the possibility of the ``label'' not being in the training set (i.e., $\mathcal{S} = \varnothing$). At inference time, for a fixed choice of guidance scale, we can combine $s^\theta_t (\bm{t}|\mathcal{S})$ with a guided vector field $u^\theta_t (\bm{t}|\mathcal{S})$ and define
\begin{align}
    \tilde{s}^\theta_\gamma (\bm{t}|\mathcal{S}) &= (1-w)s^\theta_\gamma(\bm{t}|\varnothing) + w s^\theta_\gamma (\bm{t}|\mathcal{S}) \\
    \tilde{u}^\theta_\gamma (\bm{t}|\mathcal{S}) &= (1-w)u^\theta_\gamma(\bm{t}|\varnothing) + w u^\theta_\gamma (\bm{t}|\mathcal{S})
\end{align}

Then we can sample via
\begin{align}
    \bm{t}_0 &\sim p_{init} = \mathcal{N}(0,I_d) \\
    d\bm{t}_\gamma &= \left[\tilde{u}^\theta_\gamma (\bm{t}_\gamma | \mathcal{S}) + \frac{\sigma_\gamma^2}{2} \tilde{s}^\theta_\gamma (\bm{t}_\gamma | \mathcal{S}) \right] d\gamma + \sigma_\gamma dW_\gamma \\
    \bm{\tau}^{\mathrm{SGS}} &= \bm{t}_1 \sim p_{data}(\cdot | \mathcal{S})
\end{align}
simulating the SDE from $\gamma = 0$ to 1. The goal is for $\bm{t}_1$ to adhere to the guiding variable $\mathcal{S}$.

\subsection{Generative model for memory effects}

Next, we demonstrate the applicability of generative SGS modeling for the memory term directly, improving upon traditional analytical approaches. We aim to develop a reduced model for the coefficients $r_k$ to approximate the true resolved coefficients $a_k$ for $k = 1,\dots,m$. The reduced dynamics are modeled as
\begin{align}
    \dot{r}_k (t) &= R_k(\hat{r}(t)) + \bm{\tau}_k^{\mathrm{SGS}}(\hat{r},{\xi})
\end{align}

We choose to represent the memory effects using a conditional distribution
\begin{align}
    \bm{\tau}^{\mathrm{SGS}}(\hat{r},{\xi}) \approx \int_0^t K_k (\hat{r}(t-s),s)\,ds \sim p_{data}(\cdot|\hat{r}).
\end{align}

\begin{figure}[H]
\centering

\begin{tikzpicture}[
    box/.style={
        draw,
        rectangle,
        rounded corners,
        align=center,
        minimum width=3cm,
        minimum height=2cm,
        thick,
        fill=cyan!5
    },
    arrow/.style={->, thick, shorten >=1pt},
]

\node[box] (initial) 
{$\dot{a}_k(t) = -\langle \tilde{\phi}_k, \mathcal{N}[\hat{s}] \rangle$ \\
    $a_k (0) = \varepsilon_k$ \\
    $1 \leq k \leq p$};

\node[box, below=3em of initial] (decompose) {
$a = (\hat{a},\tilde{a})$ \\
$\hat{a} = \{a\}_{k=1}^m$ \\
$\tilde{a} = \{a\}_{k=m+1}^p$};

\node[box, below=3em of decompose] (mz) {
    $\dot{r}_k(t) = R_k(\hat{r}(t)) + \bm{\tau}_k^{\mathrm{SGS}}(\hat{r},\xi)$ \\
    $R_k(\hat{r}) = -\left\langle \tilde{\phi}_k, \mathcal{N}\left[\sum_{i=1}^m \hat{r}_i \tilde{\phi}_i\right]\right\rangle$ \\
    $\bm{\tau}^{\mathrm{SGS}}(\hat{r},\xi) \sim p_{data}(\cdot | \hat{r})$
};

\draw[arrow] (initial) -- (decompose) node[midway, fill=white] {decompose};
\draw[arrow] (decompose) -- (mz) node[midway, fill=white] {SGS};

\end{tikzpicture}
\caption{Technical approach of generative SGS modeling for system of ODEs for the physics-informed expansion coefficients.}
\label{fig:mz_sgs_schematic}

\end{figure}
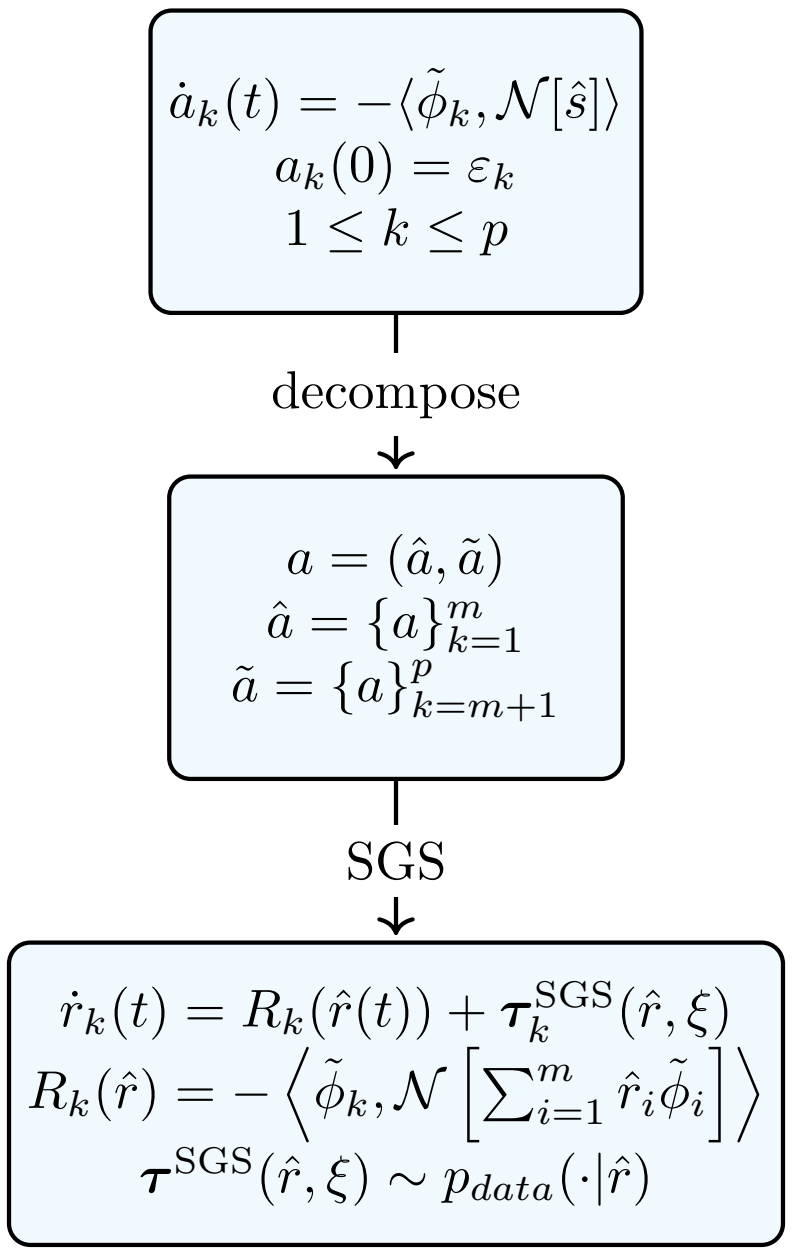

This approach differs from the previous formulation in Eq. \eqref{eq:tau_integral}, where we explicitly construct intrusive memory kernels and then evaluate the memory integral. We use diffusion modeling to learn the conditional distribution $p_{data}(\cdot|\hat{r})$, i.e., the resolved modes are conditioned on the complete set of resolved variables. This conditioning adheres with the Markovian dynamics, which are dependent on the set of resolved variables. Then, the reduced approximation $\hat{s}$ is compared to the filtered full-order field $\overline{s}$ (i.e., using only $m$ modes in the reconstruction).
\begin{align}
    \hat{s} &= \sum_{i=1}^m r_i (t) \tilde{\phi}_i (x) \\
    \bar{s} &= \sum_{i=1}^m {a}_i (t) \tilde{\phi}_i (x)
\end{align}
The technical approach is summarized in the schematic in Fig.~\ref{fig:mz_sgs_schematic}.

\begin{figure}[H]
    \centering
    \includegraphics[width=0.8\linewidth]{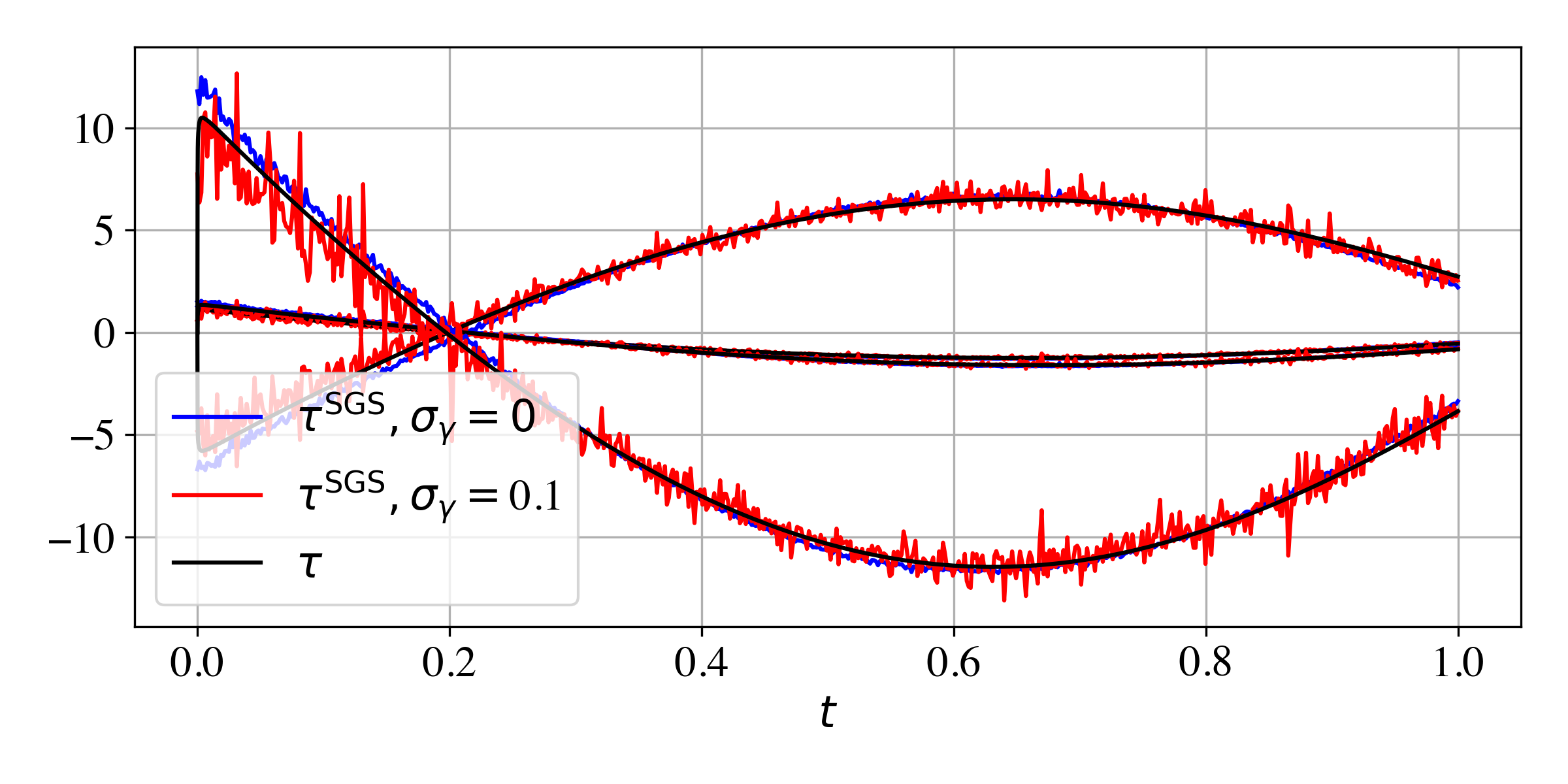}
    \caption{Generated closure $\tau_k^{SGS}$ with flow and score matching.}
    \label{prom-gen-sgs}
\end{figure}

For the reduced model, we take $m = 4$, compared to the full-order system with $p = 67$. The generative SGS terms conditioned on the test data compared to the true SGS terms are plotted in Figure \ref{prom-gen-sgs}. With the incorporation of stochastic dynamics through score matching and any $\sigma_\gamma > 0$, the resulting predictions are expected to exhibit additional variability relative to flow matching alone. For these experiments, we empirically select $\sigma_\gamma = 0.1$, which controls the magnitude of the injected diffusion and therefore the stochasticity of the predictions.

For deploying the reduced model, we use a Runge-Kutta (RK) time integration scheme at each time $t_n$, with a factor of 20 fewer timesteps compared to the full-order simulation. We sample $\bm{\tau}^{\mathrm{SGS}} (\hat{r}^n,\xi) \sim p_{data}(\cdot|\hat{r}^n)$ by simulating the SDE and then integrate forward as
\begin{align}
    \hat{r}^{n+1} = \hat{r}^{n} + \frac{h}{6}(k_1 + 2k_2 + 2k_3 + k_4) + h \bm{\tau}^{\mathrm{SGS}} (\hat{r}^n,\xi)
\end{align}
where
\begin{align}
    k_1 &= \bm{R}\left(\hat{r}^{n} \right) \\
    k_2 &= \bm{R}\left(\hat{r}^{n} + h \frac{k_1}{2}\right) \\
    k_3 &= \bm{R}\left(\hat{r}^{n} + h \frac{k_2}{2}\right) \\
    k_4 &= \bm{R}\left(\hat{r}^{n} + h k_3\right)
\end{align}

\begin{figure}[H]
    \centering
    \includegraphics[width=0.8\linewidth]{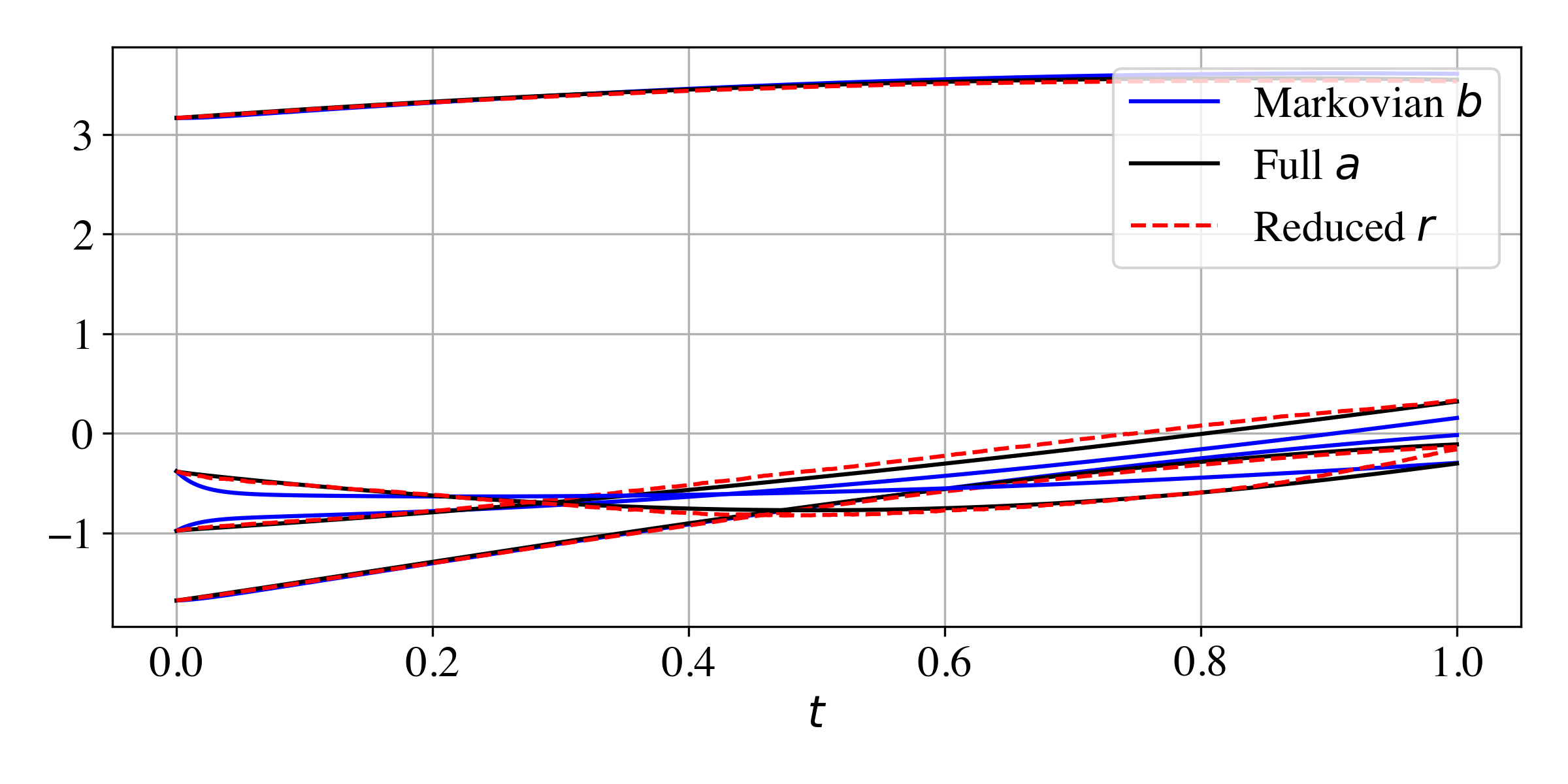}
    \caption{Full-order $a$, generated reduced $r$, and Markovian $b$ coefficients.}
    \label{prom-gen-coeffs}
\end{figure}

\begin{figure}[H]
    \centering
    \includegraphics[width=\linewidth]{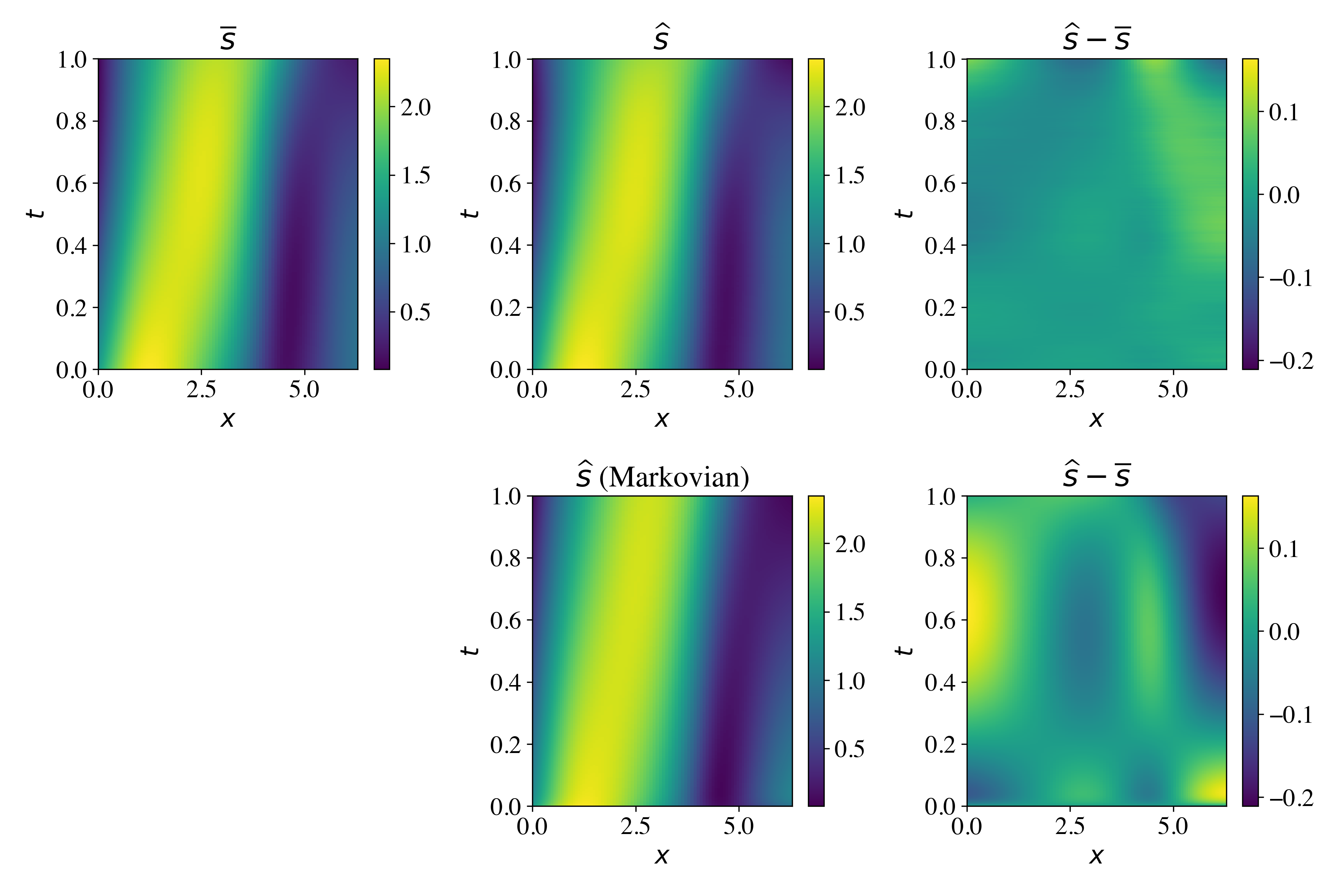}
    \caption{True filtered field $\overline{s}$ compared to Markovian and reduced approximations $\hat{s}$ with absolute error fields for initial condition $s_0 = e^{\sin(x)}$.}
    \label{prom-gen-field}
\end{figure}

For fixed function evaluation cost, we use 20 diffusion timesteps for the SDE integration at each physical timestep. Figure \ref{prom-gen-coeffs} plots the true full-order coefficients $a$, Markovian coefficients $b$, and reduced coefficients $r$. The reduced dynamics better match the true coefficients than the Markovian dynamics alone. Figure \ref{prom-gen-field} shows the reduced and Markovian fields $\hat{s}$ with absolute errors compared to the true filtered field. We see smaller absolute error magnitudes using the generative SGS model compared to the Markovian dynamics alone.

\begin{figure}[H]
    \centering
    \includegraphics[width=0.4\linewidth]{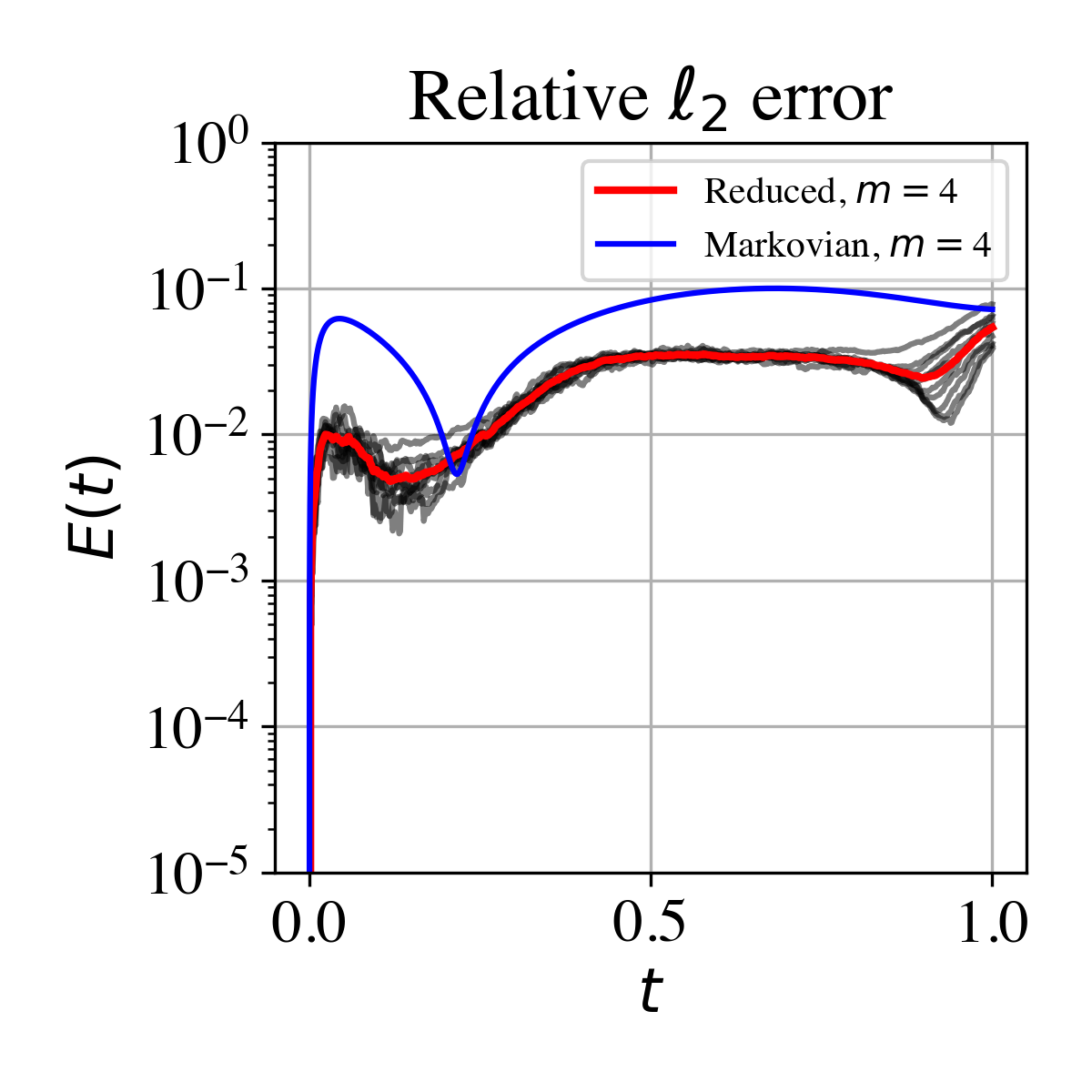}
    \caption{Relative $\ell_2$ error in time (black: samples, red: average).}
    \label{prom-gen-L2}
\end{figure}

Figure \ref{prom-gen-L2} plots the relative $\ell_2$ error of the reduced and Markovian approximations compared to the true filtered field. For this realization, we see a decreased error across the temporal domain using the generative SGS model. Table \ref{tab:viscous-burgers-prom-error} shows the average relative $\ell_2$ error, less than 2.5\% for the reduced model compared to almost 7\% for the Markovian model across 10 samples.

\begin{table}[H]
    \centering
    \begin{tabular}{|c|c|}\hline
         & Average relative $\ell_2$ error \\ \hline\hline
    Markovian & 6.6\% \\ \hline
    Reduced & 2.4\% \\ \hline
    \end{tabular}
    \caption{Average errors for generative closure model for viscous Burgers'.}
    \label{tab:viscous-burgers-prom-error}
\end{table}

\section{Conclusions}

In this work, we built reduced models by applying the MZ formalism to the system of ODEs for the physics-informed expansion coefficients. We implemented the FDT to solve for memory kernels applicable to evolving a collection of initial conditions and demonstrate lower errors on average for evolution of an out-of-distribution initial condition. These results show that higher-order approximations for the memory kernels become relevant and may be required to continue accounting for accurate memory for longer times. So, future work might include developing different representations for the memory kernels themselves, through generative approaches and higher-order projections. Finally, we used generative modeling to learn a conditional distribution to represent the memory effects on the resolved variables. We see robust predictions for evolution of a particular initial condition, with lower errors than evolving the Markovian dynamics alone. We see benefits in representing the memory itself stochastically. Future work should include accounting for memory by conditioning on previous states for some memory horizon length, as well as training across a collection of initial conditions.

\section*{Acknowledgements}

This project was completed with support from the U.S. Department of Energy, Advanced Scientific Computing Research program, under the Scalable, Efficient and Accelerated Causal Reasoning Operators, Graphs and Spikes for Earth and Embedded Systems (SEA-CROGS) project (Project No. 80278). The computational work was performed using PNNL Institutional Computing at Pacific Northwest National Laboratory. Pacific Northwest National Laboratory (PNNL) is a multi-program national laboratory operated for the U.S. Department of Energy (DOE) by Battelle Memorial Institute under Contract No. DE-AC05-76RL01830.\\

\noindent
This material is based upon work supported by the U.S. Department of Energy, Office of Science, Office of Advanced Scientific Computing Research, Department of Energy Computational Science Graduate Fellowship under Award Number DE-SC0023112. This report was prepared as an account of work sponsored by an agency of the United States Government. Neither the United States Government nor any agency thereof, nor any of their employees, makes any warranty, express or implied, or assumes any legal liability or responsibility for the accuracy, completeness, or usefulness of any information, apparatus, product, or process disclosed, or represents that its use would not infringe privately owned rights. Reference herein to any specific commercial product, process, or service by trade name, trademark, manufacturer, or otherwise does not necessarily constitute or imply its endorsement, recommendation, or favoring by the United States Government or any agency thereof. The views and opinions of authors expressed herein do not necessarily state or reflect those of the United States Government or any agency thereof.

\appendix

\section{Weak formulation derivation}
\label{sec:weak-formulation}

Rewrite the viscous Burgers' equation as
\begin{align}\label{eq:vB}
    \frac{\partial s}{\partial t} + \frac{\partial}{\partial x}\left(\frac{s^2}{2}\right) - \nu \frac{\partial^2 s}{\partial x^2} = 0
\end{align}

Set $f(s) = s^2 / 2$ and substitute $f(s)$ into Eq. \eqref{eq:vB} to write
\begin{align}
    \frac{\partial s}{\partial t} + \frac{\partial f}{\partial x} - \nu \frac{\partial^2 s}{\partial x^2} = 0
\end{align}

Rewriting as a system of first-order equations gives
\begin{align}
    s_t + f_x - \nu q_x &= 0 \\
    q - s_x &= 0
\end{align}

For the weak formulation, we can write
\begin{align}
    \langle \phi_k, s_t + f_x - \nu q_x \rangle &= 0 \label{eq:w1}\\
    \langle \phi_k, q - s_x \rangle &= 0 \label{eq:w2}
\end{align}

Integrating by parts on Eq. \eqref{eq:w1} gives
\begin{align}
    \langle \phi_k, s_t \rangle &= - \langle \phi_k, f_x \rangle + \langle \phi_k, \nu q_x \rangle \\
    \langle \phi_k, s_t \rangle &= -\bigg[\phi_k (\eta) f(\eta)\bigg]\bigg|_{\partial \Omega} + \int_{\Omega} \left(\frac{\partial}{\partial x} \phi_k (x) f(x)\right) \,dx \\ \nonumber &\hspace{5em} + \nu \bigg[\phi_k (\eta) q(\eta)\bigg]\bigg|_{\partial \Omega} - \nu \int_{\Omega} \left(\frac{\partial}{\partial x} \phi_k (x) q(x)\right) \,dx
\end{align}

Writing the flux terms as surface integrals gives
\begin{align}
    \langle \phi_k, s_t \rangle = \left\langle \phi_k', f \right\rangle - \int_{\partial \Omega} \hat{\mathbf{n}} \cdot f^*(s_R,s_L) \phi_k (\eta) \,d \eta - \nu \left\langle \phi_k', q \right\rangle + \nu \int_{\partial \Omega} \hat{\mathbf{n}} \cdot q(\eta) \phi_k (\eta) \,d \eta \label{eq:w1exp}
\end{align}
where $\hat{\mathbf{n}}$ is the outward-facing normal, $f^*$ represents the numerical flux term \cite{meuris_machine-learning-based_2023}, and the subscripts and $R$ and $L$ refer to values of the functions at the right and left boundaries. The numerical flux $f^* (c^-, c^+)$ used is defined by
\begin{align}
    f^* (c^-, c^+) &= f\left(\frac{c^- + |c^-|}{2} - \frac{c^+ - |c^+|}{2}\right) \\
    &= \begin{cases}
        f(c^-) \quad\text{if } c^+ \geq 0 \\
        f(c^+) \quad\text{if } c^- \leq 0
    \end{cases}
\end{align}

Integrating by parts on Eq. \eqref{eq:w2} gives
\begin{align}
    \langle \phi_k, q \rangle &= \langle \phi_k, s_x \rangle \\
    \langle \phi_k, q \rangle &= \bigg[\phi_k(\eta) s(\eta)\bigg]_{\partial \Omega} - \int_{\Omega} \left(\frac{\partial}{\partial x} \phi_k (x) s(x)\right) \,dx \\
    \langle \phi_k, q \rangle &= -\left\langle \phi_k', s \right\rangle + \int_{\partial \Omega} \hat{\mathbf{n}} \cdot s(\eta) \phi_k (\eta) \,d \eta \label{eq:w2exp}
\end{align}

Discretize the approximations using linear expansions in terms of the basis functions
\begin{align}
    \hat{s}_p(t,x) = \sum_{i=1}^p a_i (t) \phi_i (x), \quad
    \hat{q}_p(t,x) = \sum_{j=1}^p b_j (t) \phi_j (x).
\end{align}

Substituting the summations into Eq. \eqref{eq:w1exp} to write
\begin{multline}
    \left\langle \phi_k, \frac{\partial}{\partial t} \left( \sum_{i=1}^p a_i (t) \phi_i \right) \right\rangle = \left\langle \phi_k', f \right\rangle - \int_{\partial \Omega} \hat{\mathbf{n}} \cdot f^*(\hat{s}_{p,R},\hat{s}_{p,L}) \phi_k (\eta) \,d \eta \\ - \nu \left\langle  \phi_k', \sum_{j=1}^p b_j (t) \phi_j \right\rangle + \nu \int_{\partial \Omega} \hat{\mathbf{n}} \cdot \left( \sum_{j=1}^p b_j(t) \phi_j (x_{out}) \right) \phi_k (\eta) \,d \eta,
\end{multline}
where $x_{out} = 2\pi$ and $x_{in} = 0$. Expanding the surface integrals gives
\begin{multline}
    \left\langle \phi_k, \sum_{i=1}^p \frac{d a_i(t)}{dt} \phi_i \right\rangle = \left\langle  \phi_k', f \right\rangle - f^* ( \hat{s}_p(t, x_{out}), \hat{s}_p(t, x_{in})) (\phi_k (x_{out}) - \phi_k (x_{in})) \\ - \nu \left\langle \phi_k', \sum_{j=1}^p b_j (t) \phi_j \right\rangle + \nu \bigg[ \phi_k (x_{out}) - \phi_k (x_{in}) \bigg] \sum_{j=1}^p b_j (t) \phi_j (x_{out})
\end{multline}

As a result, the expansion coefficients for Eq. \eqref{eq:w1} are given by
\begin{multline}
    \frac{d a_k (t)}{dt} = \left\langle  \phi_k', f\right\rangle - f^* ( \hat{s}_p(t, x_{out}), \hat{s}_p(t, x_{in})) (\phi_k (x_{out}) - \phi_k (x_{in})) \\ + \nu \sum_{j=1}^p \bigg[ -\left\langle  \phi_k', \phi_j \right\rangle + \big[ \phi_k (x_{out}) - \phi_k (x_{in}) \big] \phi_j (x_{out}) \bigg] b_j (t) \label{eq:exp1}
\end{multline}
for $1 \leq k \leq p$. Substituting the summations into Eq. \eqref{eq:w2exp}
\begin{align}
    \left\langle \phi_k, \sum_{j=1}^p b_j (t) \phi_j \right\rangle = -\left\langle \phi_k', \sum_{i=1}^p a_i (t) \phi_i \right\rangle + \int_{\partial \Omega} \hat{\mathbf{n}} \cdot \left( \sum_{i=1}^p a_i (t) \phi_i (x_{in}) \right) \phi_k (\eta) \,d \eta
\end{align}

Expanding the surface integrals yields
\begin{align}
    b_k (t) = -\left\langle \phi_k', \sum_{i=1}^p a_i (t) \phi_i \right\rangle + \bigg[ \phi_k (x_{out}) - \phi_k (x_{in}) \bigg] \sum_{i=1}^p a_i (t) \phi_i (x_{in})
\end{align}

So the expansion coefficients for Eq. \eqref{eq:w2} are given by
\begin{align}
    b_k (t) = \sum_{i=1}^p \bigg[ -\left\langle \phi_k', \phi_i \right\rangle + \big[ \phi_k (x_{out}) - \phi_k (x_{in}) \big] \phi_i (x_{in}) \bigg] a_i (t) \label{eq:exp2}
\end{align}
for $1 \leq k \leq p$. Substituting Eq. \eqref{eq:exp2} into Eq. \eqref{eq:exp1} yields a system of ODEs for the expansion coefficients
\begin{multline}
    \frac{d a_k (t)}{dt} = \left\langle \phi_k', f\right\rangle - f^* ( \hat{s}_p(t, x_{out}), \hat{s}_p(t, x_{in})) (\phi_k (x_{out}) - \phi_k (x_{in})) \\ + \nu \sum_{j=1}^p \bigg[ -\left\langle \phi_k', \phi_j \right\rangle + \big[ \phi_k (x_{out}) - \phi_k (x_{in}) \big] \phi_j (x_{out}) \bigg] \cdot \\
    \sum_{i=1}^p \bigg[ -\left\langle \frac{d}{d x} \phi_j, \phi_i \right\rangle + \big[ \phi_j (x_{out}) - \phi_j (x_{in}) \big] \phi_i (x_{in}) \bigg] a_i (t)
\end{multline}
with the initial condition
\begin{align}
    a_k (0) = \langle \phi_k, s_0 \rangle
\end{align}
for $1 \leq k \leq p$. We use a quadrature rule to calculate the inner products.

Let $\mathcal{A}_k[\hat{s}_p]$ represent the nonlinear terms, $\mathcal{D}$ represent the viscous terms, and $\mathcal{Q}$ represent the derivative terms, allowing us to write
\begin{align}
    \frac{d a_k (t)}{dt} = \mathcal{A}_{k}[\hat{s}_p] + \nu \sum_{j=1}^p \mathcal{D}_{kj} \sum_{i=1}^p \mathcal{Q}_{ji} a_i (t)
\end{align}
where
\begin{align}
    \mathcal{A}_{k} &= \langle \phi'_k, f \rangle - f^* ( \hat{s}_p(t, x_{out}), \hat{s}_p(t,x_{in}))(\phi_k (x_{out}) - \phi_k(x_{in})) \\
    \mathcal{D}_{kj} &= -\langle \phi'_k, \phi_j \rangle + \left[\phi_k (x_{out}) - \phi_k (x_{in})\right]\phi_j(x_{out}) \\
    \mathcal{Q}_{ji} &= -\langle \phi'_j, \phi_i \rangle + \left[\phi_j (x_{out}) - \phi_j (x_{in})\right] \phi_i (x_{in})
\end{align}

The Markovian term is given by
\begin{align}
    R_k ({a}(\varepsilon,t)) = -\left\langle \tilde{\phi}_k, \mathcal{N}\left[\sum_{i=1}^m {a}_i(\varepsilon,t) \tilde{\phi}_i \right] \right\rangle
\end{align}
with the boundary conditions incorporated as described above. Using $\mathcal{A}, \mathcal{D},$ and $\mathcal{Q}$, we can write
\begin{align}
    \frac{d a_k (t)}{dt} = \mathcal{A}_{k}[\hat{s}_m] + \nu \sum_{j=1}^p \mathcal{D}_{jk} \sum_{i=1}^m \mathcal{Q}_{ij} a_i (t)
\end{align}
for $1 \leq k \leq m$.

\section{Fluctuations and kernels for MZ-reduced model}
\label{sec:mz-kernels}

From Appendix \ref{sec:weak-formulation}, we have the system of ODEs for the expansion coefficients
\begin{align}
    \frac{d a_k (t)}{dt} &= R_k(a(t))
\end{align}
\begin{multline}
    \frac{d a_k (t)}{dt} = \left\langle \frac{d}{d x} \phi_k, F\right\rangle - F^* ( s(t, 2\pi), s(t, 0)) (\phi_k (2\pi) - \phi_k (0)) \\ + \nu \sum_{j=1}^p \bigg[ -\left\langle \frac{d}{d x} \phi_k, \phi_j \right\rangle + \big[ \phi_k (2\pi) - \phi_k (0) \big] \phi_j (2\pi) \bigg] \cdot \\
    \sum_{i=1}^p \bigg[ -\left\langle \frac{d}{d x} \phi_j, \phi_i \right\rangle + \big[ \phi_j (2\pi) - \phi_j (0) \big] \phi_i (0) \bigg] a_i (t)
\end{multline}
with the initial condition
\begin{align}
    \varepsilon_k = a_k (0) = \langle \phi_k, s_0 \rangle
\end{align}
where
\begin{align}
    F(v) &= \frac{1}{2}v^2 \\
    F^*(c^-,c^+) &= \frac{1}{2}(\max(c^-,0) - \min(c^+,0))^2
\end{align}

For $m < p$, we want to build a ROM based on only $\hat{a} = \{a_k\}_{k=1}^m$. With initial condition $a(0) = \{\{\hat{a}_k\}_{k=1}^m, \bm{0}_{k=m+1}^p\}$. The MZ identity is
\begin{align}
    a'_k(t) = R_k((\hat{a},0)) + \sum_{j \in I}\int_0^t h_j(\hat{a}(t-s)) \kappa_{jk}(s) ds, \quad 1 \leq k \leq m
\end{align}
where $\{h_j(\hat{\varepsilon})\}_{j \in I}$ is a basis for functions of the resolved variables and $\{\kappa_{jk}(t)\}_{j \in I, 1 \leq k \leq m}$ are the memory kernels. To calculate the memory kernels \cite{chorin_optimal_2002}, we solve
\begin{align}
    \kappa_{lk}(t) = f_{lk}(t) - \sum_{j \in I} \int_0^t g_{lj}(t-s) \kappa_{jk}(s) ds
\end{align}
where
\begin{align}
    f_{lk}(t) &= \langle h_l, e^{t\mathcal{L}} \mathcal{L} \mathcal{Q} \mathcal{L} \varepsilon_k \rangle_{\mathcal{C}} \\
    g_{lj}(t) &= \langle h_l, e^{t\mathcal{L}} \mathcal{L} h_j \rangle_{\mathcal{C}}
\end{align}
where
\begin{align}
    \mathcal{L} = \sum_{j=1}^p R_j(\varepsilon) \frac{\partial}{\partial \varepsilon_j}
\end{align}

Denoting by $\{H_l(z)\}_{l \geq 0}$ the one-dimensional Hermite polynomials corresponding to weight functions $e^{-z^2/2}$, we take $\{h_j(\hat{\varepsilon})\}_{j \in I}$ to be their tensor products. We consider the case of linear polynomials $d=1$, so have $|I| = \begin{pmatrix} m+d \\ d \end{pmatrix} = m+1$. Write out the terms in $g_{lj}$
\begin{align}
    [\mathcal{L}h_j](\varepsilon) &= \frac{1}{\sigma_j} R_j (\varepsilon) \\
    e^{t\mathcal{L}} \mathcal{L} h_j &= \frac{1}{\sigma_j} R_j (a(\varepsilon,t)) \\
    g_{lj}(t) &= \frac{1}{\sigma_j} \langle h_l, R_j(a((\cdot,0),t)) \rangle_{\mathcal{C}}
\end{align}

Write out the terms in $f_{lk}$
\begin{align}
    Q \mathcal{L} \varepsilon_k &= R_k (\varepsilon) - R_k((\hat{\varepsilon},0)) \\
    \mathcal{L} Q \mathcal{L} \varepsilon_k &= \sum_{k=1}^p R_l(\varepsilon) \frac{\partial}{\partial \varepsilon_l} (Q\mathcal{L} \varepsilon_k)
\end{align}

Then we have
\begin{multline}
    e^{t\mathcal{L}} \mathcal{L} Q \mathcal{L} x_k = \langle \phi'_k, F'(s(t,\cdot)) s_t(t,\cdot) - F'(s^{m}(t,\cdot))s_t^m(t,\cdot) \rangle \\+ \nu \sum_{l=1}^p \bigg[ -\left\langle \frac{d}{d x} \phi_k, \phi_l \right\rangle + \big[ \phi_k (2\pi) - \phi_k (0) \big] \phi_l (2\pi) \bigg] \cdot \\
    \sum_{i=m+1}^p \bigg[ -\left\langle \frac{d}{d x} \phi_l, \phi_i \right\rangle + \big[ \phi_l (2\pi) - \phi_l (0) \big] \phi_i (0) \bigg] R_i (a(t)) \\ - (\phi_k(2\pi) - \phi_k(0))\bigg[ F^*_1 (s(t,2\pi), s(t,0)) s_t(t,2\pi) + F_2^* (s(t,2\pi), s(t,0)) s_t (t,0) \\- F_1^* (s^m(t,2\pi), s^m(t,0)) s_t^m(t,2\pi) - F_2^* (s^m(t,2\pi), s^m(t,0)) s_t^m(t,0)\bigg]
\end{multline}
with
\begin{align}
    F'(v) &= v \\
    F_1^* (c^-, c^+) &= (\max(c^-,0) - \min(c^+,0)) \chi_{\mathbb{R}_{>0}}(c^-) \\
    F_2^* (c^-, c^+) &= -(\max(c^-,0) - \min(c^+,0)) \chi_{\mathbb{R}_{<0}}(c^+)
\end{align}
where
\begin{align}
    \chi_{\mathbb{R}_{>0}}(c^-) &= \begin{cases}
        1 \quad\text{if } c^- > 0 \\
        0 \quad\text{if } c^- \leq 0
    \end{cases}\\
    \chi_{\mathbb{R}_{<0}}(c^+) &= \begin{cases}
        1 \quad\text{if } c^+ < 0 \\
        0 \quad\text{if } c^+ \geq 0
    \end{cases}
\end{align}

Then we can find $f_{lk}$ by computing the inner product with $\{h_l\}$. To calculate the memory kernels $\bm{K}(t) = (\kappa_{lk}(t))_{l \in I, 1 \leq k \leq m}$, discretize the convolution integral using the trapezoidal rule. At each time $t_n$, we have
\begin{align}
    \bm{K}(t_n) &= \bm{f}(t_n) - \int_0^{t_n} \bm{g}(t_m - s) \bm{K}(s) \,ds \\
    &\approx \bm{f}(t_n) - \frac{\Delta t}{2} \left[ \bm{g}(t_0) \bm{K}(t_n) + \bm{g}(t_n) \bm{K}(t_0) + 2 \sum_{q=1}^{n-1} \bm{g}(t_q) \bm{K}(t_{n-q}) \right] \\
    \left( \mathbb{I} + \frac{\Delta t}{2} \bm{g}(t_0) \right) \bm{K}(t_n) &= \bm{f}(t_n) - \frac{\Delta t}{2} \left[ \bm{g}(t_n) \bm{K}(t_0) + 2 \sum_{q=1}^{n-1} \bm{g}(t_q) \bm{K}(t_{n-q}) \right] \\
    \bm{K}(t_n) &= \left( \mathbb{I} + \frac{\Delta t}{2} \bm{g}(t_0) \right)^{-1}\left(\bm{f}(t_n) - \frac{\Delta t}{2} \left[ \bm{g}(t_n) \bm{K}(t_0) + 2 \sum_{q=1}^{n-1} \bm{g}(t_q) \bm{K}(t_{n-q}) \right]\right)
\end{align}
with $\bm{K}(t_0) = \bm{f}(t_0)$.

\section{Generative model training details and hyperparameters}
\label{sec:training}

For generative models, such as diffusion models, future work will include more rigorous assessment of hyperparameter settings, including network sizes, learning rates, and training iterations, as well as inference parameters during deployment, such as guidance scales and noise constants, as in \cite{chang_design_2026, chen_opportunities_2024, yang_diffusion_2024, chung_cfg_2024, kynkaanniemi_applying_2024}. Future studies will investigate implications of parameter choices and sensitivities, as well as assessing optimal tuning strategies for training and inference hyperparameters.

\begin{table}[H]
    \centering
    \begin{tabular}{|c|c|c|}\hline
        & Flow & Score \\\hline\hline
        Network size & $3 \times 128$ & $4 \times 256$ \\\hline
        Learning rate & 0.001 & 0.0001 \\\hline
        Activation function & SiLU & SiLU \\\hline
        Training epochs & 5000 & 5000 \\\hline
        Batch size & 100 & 50 \\\hline
        Unconditional training probability & 0.1 & 0.2 \\\hline
        Guidance scale & 1.5 & 0.8 \\\hline
    \end{tabular}
    \caption{Generative closure model training parameters}
    \label{}
\end{table}

\printbibliography

\end{document}